\documentclass{amsart}
\usepackage{amsmath,amssymb,amsthm,bm,bbm}
\usepackage{amscd}
\usepackage{graphicx,enumerate}
\usepackage{multirow}
\usepackage{layout}
\usepackage[OT2,T1]{fontenc}
\DeclareSymbolFont{cyrletters}{OT2}{wncyr}{m}{n}
\DeclareMathSymbol{\Sha}{\mathalpha}{cyrletters}{"58}
\usepackage[OT2,OT1]{fontenc}
\newcommand\cyr{\renewcommand\rmdefault{wncyr}
\renewcommand\sfdefault{wncyss}
\renewcommand\encodingdefault{OT2}
\normalfont\selectfont}
\DeclareTextFontCommand{\textcyr}{\cyr}

\theoremstyle{plain}
\newtheorem{theorem}{Theorem}[section]
\newtheorem*{theorem-nn}{Theorem}

\newtheorem{proposition}[theorem]{Proposition}
\newtheorem*{proposition-nn}{Proposition}

\theoremstyle{definition}
\newtheorem{definition}[theorem]{Definition}
\newtheorem{example}[theorem]{Example}
\newtheorem{remark}[theorem]{Remark}

\newtheorem*{acknowledgments}{Acknowledgments}

\theoremstyle{remark}

\newcommand{\bZ}{\mathbbm{Z}}\newcommand{\bQ}{\mathbbm{Q}}
\newcommand{\bC}{\mathbbm{C}}
\newcommand{\bG}{\mathbbm{G}}\newcommand{\bF}{\mathbbm{F}}
\newcommand{\bA}{\mathbbm{A}}
\newcommand{\cC}{\mathcal{C}}\newcommand{\cD}{\mathcal{D}}
\newcommand{\cH}{\mathcal{H}}\newcommand{\cS}{\mathcal{S}}

\newcommand{\PSL}{{\rm PSL}}

\title{Norm one tori and Hasse norm principle, II: Degree $12$ case}
\author[A. Hoshi]{Akinari Hoshi}
\address{Department of Mathematics, Niigata University, Niigata 950-2181, Japan}
\email{hoshi@math.sc.niigata-u.ac.jp}

\author[K. Kanai]{Kazuki Kanai}
\address{Graduate School of Science and Technology, Niigata University, Niigata 950-2181, Japan}
\email{kanai@m.sc.niigata-u.ac.jp}

\author[A. Yamasaki]{Aiichi Yamasaki}
\address{Department of Mathematics, Kyoto University, Kyoto 606-8502, Japan}
\email{aiichi.yamasaki@gmail.com}

\thanks{%, flabby resolution.\\
This work was partially supported by JSPS KAKENHI Grant Numbers 
16K05059, 19K03418, 20K03511. 
Parts of the work were finished when the
first-named author and the third-named author 
were visiting the National Center for Theoretic Sciences (Taipei),
whose support is gratefully acknowledged.
}

\begin{document}
\maketitle
\begin{abstract}
Let $k$ be a field, $T$ be an algebraic $k$-torus, 
$X$ be a smooth $k$-compactification of $T$ 
and ${\rm Pic}\,\overline{X}$ be the Picard group of 
$\overline{X}=X\times_k\overline{k}$. 
Hoshi, Kanai and Yamasaki \cite{HKY} determined 
$H^1(k,{\rm Pic}\, \overline{X})$ 
for norm one tori $T=R^{(1)}_{K/k}(\bG_m)$ 
and gave 
a necessary and sufficient condition for the Hasse norm principle for 
extensions $K/k$ of number fields with $[K:k]=n\leq 15$ and $n\neq 12$. 
In this paper, we determine $64$ cases with 
$H^1(k,{\rm Pic}\, \overline{X})\neq 0$
and give a necessary and sufficient condition for the Hasse norm principle 
for $K/k$ with $[K:k]=12$. 
\end{abstract}
\tableofcontents
%
%%%%%%%%%%%%%%%%%%%%%%%%%%%%%%%%%%%%%%%%%%%%%
\section{Introduction}\label{S1}
Let $k$ be a field, 
$K/k$ be a separable field extension of degree $n$ 
and $L/k$ be the Galois closure of $K/k$ with ${\rm Gal}(L/k)\simeq G$. 
Let $S_n$ (resp. $A_n$, $D_n$, $C_n$) be the symmetric 
(resp. the alternating, the dihedral, the cyclic) group 
of degree $n$ of order $n!$ (resp. $n!/2$, $2n$, $n$). 
Let $nTm$ be the $m$-th transitive subgroup of the symmetric group 
$S_n$ of degree $n$ up to conjugacy 
(see Butler and McKay \cite{BM83}, \cite{GAP}). 
We may regard $G$ as the transitive subgroup $G=nTm\leq S_n$. 
Let $H={\rm Gal}(L/K)\leq G$ with $[G:H]=n$. 
Then we may assume that 
$H$ is the stabilizer of one of the letters in $G$, 
i.e. $L=k(\theta_1,\ldots,\theta_n)$ and $K=k(\theta_i)$ for some 
$1\leq i\leq n$.

Voskresenskii \cite{Vos67} proved that 
all the $2$-dimensional algebraic $k$-tori are $k$-rational. 
This implies that 
$H^1(k,{\rm Pic}\,\overline{X})\simeq A(T)\simeq \Sha(T)=0$ 
where $X$ is a smooth $k$-compactification of an algebraic $k$-torus $T$, 
$A(T)$ is the kernel of the weak approximation of $T$ 
and $\Sha(T)$ is the Shafarevich-Tate group of $T$ 
(see Section \ref{S3} and also Manin \cite[\S 30]{Man74}). 
Kunyavskii \cite{Kun84} showed that, among $73$ cases of 
$3$-dimensional $k$-tori $T$, 
there exist exactly $2$ cases with 
$H^1(k,{\rm Pic}\,\overline{X})\neq 0$ 
which are of special type called norm one tori. 
Hoshi, Kanai and Yamasaki \cite[Theorem 1.18]{HKY} determined 
$H^1(k,{\rm Pic}\, \overline{X})$ 
for norm one tori $T=R^{(1)}_{K/k}(\bG_m)$ 
with $[K:k]=n\leq 15$ and $n\neq 12$. 
%Note that there exist $2$ (resp. $5$, $5$, $16$, $7$, $50$, $34$, $45$, $8$, 
%$301$, $9$, $63$, $104$) transitive subgroups $nTm$ of $S_n$ up to conjugacy 
%for $n=3$ (resp. $4$, $5$, $6$, $7$, $8$, $9$, $10$, $11$, $12$, 
%$13$, $14$, $15$) 
%(see Butler and McKay \cite{BM83} for $n\leq 11$, 
%Royle \cite{Roy87} for $n=12$, 
%Butler \cite{But93} for $n=14,15$ and \cite{GAP}). 
We determine $64$ cases with 
$H^1(k,{\rm Pic}\, \overline{X})\neq 0$ when $[K:k]=12$ as follows: 
(Note that there exist exactly $301$ transitive subgroups $12Tm$ 
of $S_{12}$ up to conjugacy (see Royle \cite{Roy87} and \cite{GAP}).) 

%%%%%%%%%%%%%%%%%%%%%%%%%%%%%%%%%%%%%%%%%%%%%%%%%%%%%%%%%%%%%%%
\begin{theorem}\label{thmain1}
Let $k$ be a field, 
$K/k$ be a separable field extension of degree $12$ 
and $L/k$ be the Galois closure of $K/k$. 
Assume that $G={\rm Gal}(L/k)=12Tm$ $(1\leq m\leq 301)$ 
is a transitive subgroup of $S_{12}$ 
and $H={\rm Gal}(L/K)$ with $[G:H]=12$. 
Let $T=R^{(1)}_{K/k}(\bG_m)$ be the norm one torus of $K/k$ of dimension $11$ 
and $X$ be a smooth $k$-compactification of $T$. 
Then $H^1(k,{\rm Pic}\, \overline{X})\neq 0$ 
if and only if $G$ is given as in Table $1$. 
In particular, if $k$ is a number field and 
$L/k$ is an unramified extension, then $A(T)=0$ and 
$H^1(k,{\rm Pic}\,\overline{X})\simeq \Sha(T)$. 
\end{theorem}

In Table $1$, 
$V_4\simeq C_2\times C_2$ is the Klein four group, 
$Q_8$ is the quaternion group of order $8$, 
$\PSL_2(\bF_{11})$ is the projective special linear group of degree $2$ 
over the finite field $\bF_{11}$ of $11$ elements, and 
$S_n(m)$ and $A_n(m)$ mean that $S_n(m)\simeq S_n=mTx\leq S_m$ 
and $A_n(m)\simeq A_n=mTx\leq S_m$ respectively. 

%%%%%%%%%%%%%%%%%%%%%%%%%%%%%%%%%%%%%%%%%%%%%%%%
%
%\newpage
\begin{center}
\vspace*{1mm}
Table $1$: $H^1(k,{\rm Pic}\, \overline{X})\simeq H^1(G,[J_{G/H}]^{fl})\neq 0$ 
with $G=12Tm$ $(1\leq m\leq 301)$\vspace*{2mm}\\
\begin{tabular}{lc} 
$G$ & $H^1(k,{\rm Pic}\, \overline{X})$ 
$\simeq H^1(G,[J_{G/H}]^{fl})$\\\hline
$12T2\simeq C_6\times C_2$ & $\bZ/2\bZ$\\
$12T3\simeq D_6$ & $\bZ/2\bZ$\\
$12T4\simeq A_4(12)$ & $\bZ/2\bZ$\\
$12T7\simeq A_4(6)\times C_2$ & $\bZ/2\bZ$\\
$12T9\simeq S_4$ & $\bZ/2\bZ$\\
$12T10\simeq S_3\times V_4$ & $\bZ/2\bZ$\\
$12T16\simeq (S_3)^2$ & $\bZ/2\bZ$\\
$12T18\simeq S_3\times C_6$ & $\bZ/2\bZ$\\
$12T20\simeq A_4(4)\times C_3$ & $\bZ/2\bZ$\\
$12T31\simeq (C_4)^2\rtimes C_3$ & $\bZ/4\bZ$\\
$12T32\simeq (C_2)^4\rtimes C_3$ & $(\bZ/2\bZ)^{\oplus 2}$\\
$12T33\simeq A_5(12)$ & $\bZ/2\bZ$\\
$12T34\simeq (S_3)^2\rtimes C_2$ & $\bZ/2\bZ$\\
$12T37\simeq (S_3)^2\times C_2$ & $\bZ/2\bZ$\\
$12T40\simeq ((C_3)^2\rtimes C_4)\times C_2$ & $\bZ/2\bZ$\\
$12T43\simeq A_4(4)\times S_3$ & $\bZ/2\bZ$\\
$12T47\simeq (C_3)^2\rtimes Q_8$ & $\bZ/2\bZ$\\
$12T52\simeq (A_4\times V_4)\rtimes C_2$ & $\bZ/2\bZ$\\
$12T54\simeq (S_4\times C_2)\rtimes C_2$ & $\bZ/2\bZ$\\
$12T55\simeq ((C_4)^2\rtimes C_3)\times C_2$ & $\bZ/2\bZ$\\
$12T56\simeq ((C_2)^4\rtimes C_3)\times C_2$ & $\bZ/2\bZ$\\
$12T57\simeq ((C_4\times C_2)\rtimes C_4)\rtimes C_3$ & $\bZ/4\bZ$\\
$12T59\simeq (C_2)^3\rtimes A_4(6)$ & $\bZ/2\bZ$\\
$12T61\simeq ((C_4)^2\rtimes C_2)\rtimes C_3$ & $\bZ/2\bZ$\\
$12T64\simeq ((C_4)^2\rtimes C_3)\rtimes C_2$ & $\bZ/2\bZ$\\
$12T65\simeq ((C_4)^2\rtimes C_3)\rtimes C_2$ & $\bZ/2\bZ$\\
$12T66\simeq ((C_2)^4\rtimes C_3)\rtimes C_2$ & $\bZ/2\bZ$\\
$12T70\simeq (S_3)^2\times C_3$ & $\bZ/2\bZ$\\
$12T71\simeq (C_3)^3\rtimes V_4$ & $\bZ/2\bZ$\\
$12T74\simeq S_5(12)$ & $\bZ/2\bZ$\\
$12T75\simeq A_5(6)\times C_2$ & $\bZ/2\bZ$\\
$12T77\simeq (S_3)^2\rtimes V_4$ & $\bZ/2\bZ$\\\hline
\end{tabular}
\end{center}~\\

%\newpage
\begin{center}
\vspace*{1mm}
Table $1$ (continued): $H^1(k,{\rm Pic}\, \overline{X})\simeq H^1(G,[J_{G/H}]^{fl})\neq 0$ 
with $G=12Tm$ $(1\leq m\leq 301)$\vspace*{2mm}\\
\begin{tabular}{lc} 
$G$ & $H^1(k,{\rm Pic}\, \overline{X})$ 
$\simeq H^1(G,[J_{G/H}]^{fl})$\\\hline
$12T88\simeq (C_2)^4\rtimes A_4(6)$ & $\bZ/2\bZ$\\
$12T92\simeq (((C_4)^2\rtimes C_2)\rtimes C_3)\times C_2$ & $\bZ/2\bZ$\\
$12T93\simeq (((C_4\rtimes C_4)\times C_2)\rtimes C_2)\rtimes C_3$ & $\bZ/2\bZ$\\
$12T96\simeq (((C_4)^2\rtimes C_3)\rtimes C_2)\times C_2$ & $\bZ/2\bZ$\\
$12T97\simeq (((C_4)^2\rtimes C_3)\rtimes C_2)\times C_2$ & $\bZ/2\bZ$\\
$12T100\simeq (((C_2)^4\rtimes C_3)\rtimes C_2)\times C_2$ & $\bZ/2\bZ$\\
$12T102\simeq ((C_2)^4\rtimes C_3)\rtimes C_4$ & $\bZ/2\bZ$\\
$12T117\simeq (S_3)^3$ & $\bZ/2\bZ$\\
$12T122\simeq ((C_3)^2\rtimes Q_8)\rtimes C_3$ & $\bZ/2\bZ$\\
$12T130\simeq (C_3)^4\rtimes V_4\simeq C_3\wr V_4$ & $\bZ/2\bZ$\\
$12T132\simeq ((C_3)^3\rtimes V_4)\rtimes C_3$ & $\bZ/2\bZ$\\
$12T133\simeq (C_3)^3\rtimes A_4(4)$ & $\bZ/2\bZ$\\
$12T144\simeq (C_2)^5\rtimes A_4(6)$ & $\bZ/2\bZ$\\
$12T168\simeq (C_3)^4\rtimes (C_2)^3$ & $\bZ/2\bZ$\\
$12T171\simeq (C_3)^4\rtimes (C_4\times C_2)$ & $\bZ/2\bZ$\\
$12T172\simeq (C_3)^4\rtimes D_4$ & $\bZ/2\bZ$\\
$12T174\simeq (C_3)^4\rtimes Q_8$ & $\bZ/2\bZ$\\
$12T176\simeq ((C_3)^3\rtimes C_2)\rtimes A_4(4)$ & $\bZ/2\bZ$\\
$12T179\simeq \PSL_2(\bF_{11})$ & $\bZ/2\bZ$\\
$12T188\simeq (C_2)^6\rtimes A_4(6)\simeq C_2\wr A_4(6)$ & $\bZ/2\bZ$\\
$12T194\simeq (C_3)^4\rtimes A_4(4)\simeq C_3\wr A_4(4)$ & $\bZ/2\bZ$\\
$12T210\simeq (C_3)^4\rtimes (D_4\times C_2)$ & $\bZ/2\bZ$\\
$12T214\simeq (C_3)^4\rtimes ((C_4\times C_2)\rtimes C_2)$ & $\bZ/2\bZ$\\
$12T230\simeq (C_2)^5\rtimes A_5(6)$ & $\bZ/2\bZ$\\
$12T232\simeq ((C_3)^4\rtimes Q_8)\rtimes C_3$ & $\bZ/2\bZ$\\
$12T234\simeq ((C_3)^4\rtimes C_2)\rtimes A_4(4)$ & $\bZ/2\bZ$\\
$12T242\simeq (C_3)^4\rtimes ((C_2)^3\rtimes V_4)$ & $\bZ/2\bZ$\\
$12T246\simeq (C_3)^4\rtimes ((C_2)^3\rtimes C_4)$ & $\bZ/2\bZ$\\
$12T255\simeq (C_2)^6\rtimes A_5(6)\simeq C_2\wr A_5(6)$ & $\bZ/2\bZ$\\
$12T261\simeq (S_3)^4\rtimes V_4\simeq S_3\wr V_4$ & $\bZ/2\bZ$\\
$12T271\simeq ((C_3)^4\rtimes (C_2)^3)\rtimes A_4(4)$ & $\bZ/2\bZ$\\
$12T280\simeq (S_3)^4\rtimes A_4(4)\simeq S_3\wr A_4(4)$ & $\bZ/2\bZ$\\\hline
\end{tabular}
\end{center}~\\
%%%%%%%%%%%%%%%%%%%%%%%%%%%%%%%%%%%%%%%%%%%%%%%%
\newpage

\begin{remark}
Theorem \ref{thmain1} enables us 
to obtain the group $T(k)/R$ of $R$-equivalence classes 
over a local field $k$ via 
$T(k)/R\simeq H^1(k,{\rm Pic}\,\overline{X})\simeq 
H^1(G,[J_{G/H}]^{fl})$ for norm one tori $T=R^{(1)}_{K/k}(\bG_m)$ 
with $[K:k]=12$ 
(see Colliot-Th\'{e}l\`{e}ne and Sansuc \cite[Corollary 5, page 201]{CTS77}, 
Voskresenskii \cite[Section 17.2]{Vos98} and Hoshi, Kanai and Yamasaki \cite[Section 7, Application 1]{HKY}).
\end{remark}

Kunyavskii \cite{Kun84} gave a necessary and sufficient condition 
for the Hasse norm principle for $K/k$ with $[K:k]=4$ $(G=4Tm\ (1\leq m\leq 5))$. 
Drakokhrust and Platonov \cite{DP87} gave a necessary and sufficient condition 
for the Hasse norm principle for $K/k$ with $[K:k]=6$ $(G=6Tm\ (1\leq m\leq 16))$. 
Hoshi, Kanai and Yamasaki \cite[Theorem 1.18]{HKY} gave 
a necessary and sufficient condition for the Hasse norm principle 
for $K/k$ with $[K:k]=n\leq 15$ and $n\neq 12$ 
(see also \cite[Section 1]{HKY}). 

By Theorem \ref{thmain1} and Ono's theorem (see Theorem \ref{thOno}) 
which claims that 
the Hasse norm principle holds for $K/k$ if and only if 
$\Sha(T)=0$ where $T=R^{(1)}_{K/k}(\bG_m)$ is the norm one torus of $K/k$, 
the Hasse norm principle holds for $K/k$ with $[K:k]=12$ as in the cases 
in Theorem \ref{thmain1} except for the cases as in Table $1$. 
By using Drakokhrust and Platonov's method 
(see Section \ref{S6} for details) 
and some new useful functions of GAP \cite{GAP} provided in Section \ref{S7}, 
we will prove the following main theorem of this paper which gives 
a necessary and sufficient condition 
for the Hasse norm principle for $K/k$ with $[K:k]=12$: 
\begin{theorem}\label{thmain2}
Let $k$ be a number field, 
$K/k$ be a field extension of degree $12$ 
and $L/k$ be the Galois closure of $K/k$. 
Let $G={\rm Gal}(L/k)=12Tm$ $(1\leq m\leq 301)$ 
be a transitive subgroup of $S_{12}$, 
$H={\rm Gal}(L/K)$ with $[G:H]=12$ 
and $G_v$ be the decomposition group of $G$ at a place $v$ of $k$. 
Let $T=R^{(1)}_{K/k}(\bG_m)$ be the norm one torus of $K/k$ 
of dimension $11$ and $X$ be a smooth $k$-compactification of $T$. 
Then $A(T)\simeq \Sha(T)=0$ except for the $64$ cases in Table $1$. 
For the $64$ cases in Table $1$ except for $3$ cases 
$G=12T31$, $12T32$, $12T57$ as in Tables $2$-$3$, $2$-$4$, $2$-$5$, 
we have either {\rm (a)} 
$A(T)=0$ and $\Sha(T)\simeq\bZ/2\bZ$ 
or {\rm (b)} 
$A(T)\simeq\bZ/2\bZ$ and $\Sha(T)=0$. 
We assume that $H$ is the stabilizer of one of the letters in $G$. 
Then $\Sha(T)$ is given as in Tables $2$-$1$ to $2$-$5$. 
Moreover, for the cases in Table $2$-$1$ and Table $2$-$3$, 
$\Sha(T)$ is also given for general $H\leq G$ with $[G:H]=12$. 
\end{theorem}

In Tables $2$-$1$ to $2$-$5$, 
$M_{16}=\langle x,y\mid x^8=y^2=1, yxy^{-1}=x^5\rangle$ 
is the modular type $2$-group of order $16$, 
$Z(G)$ is the center of a group $G$, 
$D(G)$ is the derived subgroup of $G$ 
and 
$D^i(G):=D(D^{i-1}(G))$ is the $i$-th derived subgroup 
$(D^0(G):=G)$, 
${\rm Syl}_p(G)$ is a $p$-Sylow subgroup of $G$ 
where $p$ is a prime number, 
$N_G(G^\prime)$ is the normalizer of 
a subgroup $G^\prime\leq G$, 
$\Phi(G)$ is the Frattini subgroup of $G$, 
i.e. the intersection of all maximal subgroups of $G$, 
and 
${\rm Orb}_{G}(i)$ is the orbit of $1\leq i\leq n$ 
under the action of $G\leq S_n$. 

Note that a place $v$ of $k$ with non-cyclic decomposition group $G_v$ 
as in Tables $2$-$1$ to $2$-$5$ must be ramified in $L$ because 
if $v$ is unramified, then $G_v$ is cyclic. 

\newpage
%%%%%%%%%%%%%%%%%%%%%%%%%%%%%%%%%%%%%%%%%%%%%%%%%%%%%%%%%%%%%%%%%%%%%
\begin{center}
\vspace*{1mm}
Table $2$-$1$: $\Sha(T)=0$ for $T=R^{(1)}_{K/k}(\bG_m)$ 
and $G={\rm Gal}(L/k)=12Tm$ as in Table $1$\vspace*{2mm}\\
\renewcommand{\arraystretch}{1.05}
\begin{tabular}{ll} 
 & $\Sha(T)\leq \bZ/2\bZ$, and\\
$G$ & $\Sha(T)=0$ if and only if there exists a place $v$ of $k$ such that\\\hline\vspace*{-4mm}\\
$12T2\simeq C_6\times C_2$ & \multirow{16}{*}{
\begin{minipage}{12cm}
$V_4\leq G_v$ {\rm (Tate \cite{Tat67} for $12T2$, $12T3$, $12T4$)}
\end{minipage}}\\
$12T3\simeq D_6$ & \\
$12T4\simeq A_4(12)$ & \\
$12T16\simeq (S_3)^2$ & \\
$12T18\simeq S_3\times C_6$ & \\
\end{tabular}
\end{center}~\\
%%%%%%%

%\newpage
%%%%%%%%%%%%%%%%%%%%%%%%%%%%%%%%%%%%%%%%%%%%%%%%%%%%%%%%%%%%%%%%%%%%%
\begin{center}
\vspace*{1mm}
Table $2$-$1$ (continued): $\Sha(T)=0$ for $T=R^{(1)}_{K/k}(\bG_m)$ 
and $G={\rm Gal}(L/k)=12Tm$ as in Table $1$\vspace*{2mm}\\
\renewcommand{\arraystretch}{1.05}
\begin{tabular}{ll} 
 & $\Sha(T)\leq \bZ/2\bZ$, and\\
$G$ & $\Sha(T)=0$ if and only if there exists a place $v$ of $k$ such that\\\hline\vspace*{-4mm}\\
%%%%%%%%%%%%%%%%%%%%%%%%%%%%%%%%%%
$12T64\simeq ((C_4)^2\rtimes C_3)\rtimes C_2$ & \multirow{4}{*}{
\begin{minipage}{12cm}
{\rm (i)} $D_4\leq G_v$ where 
$\Phi(D^2(G))\leq D_4$ 
for the unique characteristic subgroup $\Phi(D^2(G))\simeq V_4\lhd G$ 
of order $4$ 
$($equivalently, $|{\rm Orb}_{D_4}(i)|=4$ for $8$ of $1\leq i\leq 12$ and 
$|{\rm Orb}_{D_4}(i)|=2$ for $4$ of $1\leq i\leq 12$$)$, or\\ 
{\rm (ii)} $(C_4)^2\leq G_v$
\end{minipage}}\\
 & \\
 & \\
 & \vspace*{1mm}\\\hline\vspace*{-4mm}\\
\end{tabular}
\end{center}~\\
%%%%%%%

%\newpage
%%%%%%%%%%%%%%%%%%%%%%%%%%%%%%%%%%%%%%%%%%%%%%%%%%%%%%%%%%%%%%%%%%%%%
\begin{center}
\vspace*{1mm}
Table $2$-$2$: $\Sha(T)=0$ for $T=R^{(1)}_{K/k}(\bG_m)$ 
and $G={\rm Gal}(L/k)=12Tm$ as in Table $1$\vspace*{2mm}\\
\renewcommand{\arraystretch}{1.05}
\begin{tabular}{ll} 
 & $\Sha(T)\leq \bZ/2\bZ$, and\\
$G$ & $\Sha(T)=0$ if and only if there exists a place $v$ of $k$ such that\\\hline\vspace*{-4mm}\\
%%%%%%%%%%%%%%%%%%%%%%%%%%%%%%%%%%
$12T10\simeq S_3\times V_4$ & \multirow{2}{*}{
\begin{minipage}{12cm}
$V_4\leq G_v$ where 
$|{\rm Orb}_{V_4}(i)|=4$ for any $1\leq i\leq 12$
\end{minipage}}\\
$12T37\simeq (S_3)^2\times C_2$ & \vspace*{1mm}\\\hline\vspace*{-4mm}\\
%%%%%%%%%%%%%%%%%%%%%%%%%%%%%%%%%%
$12T54\simeq (S_4\times C_2)\rtimes C_2$ & \multirow{6}{*}{
\begin{minipage}{12cm}
{\rm (i)} $V_4\leq G_v$ where 
$|{\rm Orb}_{V_4}(i)|=2$ for any $1\leq i\leq 12$ 
and $V_4\cap D(G)\simeq C_2$ with $D(G)\simeq A_4\times C_2$,\\ 
{\rm (ii)} $C_4\times C_2\leq G_v$ where 
$|{\rm Orb}_{C_4\times C_2}(i)|=4$ for $8$ of $1\leq i\leq 12$ and 
$|{\rm Orb}_{C_4\times C_2}(i)|=2$ for $4$ of $1\leq i\leq 12$, or\\ 
{\rm (iii)} $D_4\leq G_v$ where 
$|{\rm Orb}_{D_4}(i)|=4$ for $8$ of $1\leq i\leq 12$, 
$|{\rm Orb}_{D_4}(i)|=2$ for $4$ of $1\leq i\leq 12$, 
and $D^2(G)\leq D_4$ 
with $D^2(G)\simeq V_4\lhd G$
\end{minipage}}\\
 & \\
\end{tabular}
\end{center}~\\
%%%%%%%

%\newpage
%%%%%%%%%%%%%%%%%%%%%%%%%%%%%%%%%%%%%%%%%%%%%%%%%%%%%%%%%%%%%%%%%%%%%
\begin{center}
\vspace*{1mm}
Table $2$-$2$ (continued): $\Sha(T)=0$ for $T=R^{(1)}_{K/k}(\bG_m)$ 
and $G={\rm Gal}(L/k)=12Tm$ as in Table $1$\vspace*{2mm}\\
\renewcommand{\arraystretch}{1.05}
\begin{tabular}{ll} 
 & $\Sha(T)\leq \bZ/2\bZ$, and\\
$G$ & $\Sha(T)=0$ if and only if there exists a place $v$ of $k$ such that\\\hline\vspace*{-4mm}\\
%%%%%%%%%%%%%%%%%%%%%%%%%%%%%%%%%%
$12T102\simeq ((C_2)^4\rtimes C_3)\rtimes C_4$ & \multirow{8}{*}{
\begin{minipage}{12cm} 
{\rm (i)} $V_4\leq G_v$ where 
$|{\rm Orb}_{V_4}(i)|=2$ for any $1\leq i\leq 12$ and 
{\rm (i-1)}
$V_4\cap D(G)\simeq C_2$ with $D(G)\simeq (C_2)^4\rtimes C_3\lhd G$, 
$V_4\cap Z(G)=1$ with $Z(G)\simeq C_2\lhd G$ and 
any non-zero element of $V_4$ is a product of $4$ $($different$)$ 
transpositions 
or, 
{\rm (i-2)} 
$V_4\leq D(G)$, $V_4\cap Z(G)=1$ and 
$V_4$ is not a normal subgroup of $G$,\\ 
{\rm (ii)} $C_4\times C_2\leq G_v$ where 
$|{\rm Orb}_{C_4\times C_2}(i)|=4$ for $8$ of $1\leq i\leq 12$ and 
$|{\rm Orb}_{C_4\times C_2}(i)|=2$ for $4$ of $1\leq i\leq 12$, or\\ 
{\rm (iii)} $G_{16,3}:=(C_4\times C_2)\rtimes C_2\leq G_v$ % G(16,3) 
where
$|{\rm Orb}_{G_{16,3}}(i)|=8$ for $8$ of $1\leq i\leq 12$ and 
$|{\rm Orb}_{G_{16,3}}(i)|=2$ for $4$ of $1\leq i\leq 12$
\end{minipage}}\\
 & \\
 & \\
 & \\
 & \\
\end{tabular}
\end{center}~\\

\newpage
%%%%%%%%%%%%%%%%%%%%%%%%%%%%%%%%%%%%%%%%%%%%%%%%%%%%%%%%%%%%%%%%%%%%%
\begin{center}
\vspace*{1mm}
Table $2$-$2$ (continued): $\Sha(T)=0$ for $T=R^{(1)}_{K/k}(\bG_m)$ 
and $G={\rm Gal}(L/k)=12Tm$ as in Table $1$\vspace*{2mm}\\
\renewcommand{\arraystretch}{1.05}
\begin{tabular}{ll} 
 & $\Sha(T)\leq \bZ/2\bZ$, and\\
$G$ & $\Sha(T)=0$ if and only if there exists a place $v$ of $k$ such that\\\hline\vspace*{-4mm}\\
%%%%%%%%%%%%%%%%%%%%%%%%%%%%%%%%%%
$12T214$ & \multirow{5}{*}{
\begin{minipage}{12cm} 
{\rm (i)} $C_4\times C_2\leq G_v$ where 
$|{\rm Orb}_{C_4\times C_2}(i)|=8$ for $8$ of $1\leq i\leq 12$ and 
$|{\rm Orb}_{C_4\times C_2}(i)|=4$ for $4$ of $1\leq i\leq 12$,\\ 
{\rm (ii)} $D_4\leq G_v$ where 
$|{\rm Orb}_{D_4}(i)|=8$ for $8$ of $1\leq i\leq 12$ and 
$|{\rm Orb}_{D_4}(i)|=4$ for $4$ of $1\leq i\leq 12$, or\\
{\rm (iii)} $Q_8\leq G_v$ 
\end{minipage}}\\
$\simeq (C_3)^4\rtimes ((C_4\times C_2)\rtimes C_2)$ & \\
 & \\
 & \\
 & \vspace*{1mm}\\
\end{tabular}
\end{center}~\\

\newpage
%%%%%%%%%%%%%%%%%%%%%%%%%%%%%%%%%%%%%%%%%%%%%%%%%%%%%%%%%%%%%%%%%%%%%
\begin{center}
\vspace*{1mm}
Table $2$-$3$: $\Sha(T)$ for $T=R^{(1)}_{K/k}(\bG_m)$ 
and $G={\rm Gal}(L/k)=12T31$ as in Table $1$\vspace*{2mm}\\
\renewcommand{\arraystretch}{1.05}
\begin{tabular}{ll} 
$G$ & $\Sha(T)\leq \bZ/4\bZ$\\\hline\vspace*{-4mm}\\
%%%%%%%%%%%%%%%%%%%%%%%%%%%%%%%%%%
$12T31\simeq (C_4)^2\rtimes C_3$\hspace*{5mm} & \multirow{4}{*}{
\begin{minipage}{12cm} 
{\rm (I)} $\Sha(T)=0$ if and only if 
there exists a place $v$ of $k$ such that 
$(C_4)^2\leq G_v$;\\
{\rm (II)} $\Sha(T)\simeq\bZ/2\bZ$ if and only if 
there exists a place $v$ of $k$ such that $V_4\leq G_v$ 
and {\rm (I)} does not hold
\end{minipage}}\\
 & \\
 & \\
 & \vspace*{1mm}\\\hline
%%%%%%%%%%%%%%%%%%%%%%%%%%%%%%%%%%
\end{tabular}
\end{center}~\\
\vspace*{5mm}

\begin{center}
\vspace*{1mm}
Table $2$-$4$: $\Sha(T)$ for $T=R^{(1)}_{K/k}(\bG_m)$ 
and $G={\rm Gal}(L/k)=12T57$ as in Table $1$\vspace*{2mm}\\
\renewcommand{\arraystretch}{1.05}
\begin{tabular}{ll} 
$G$ & $\Sha(T)\leq \bZ/4\bZ$\\\hline\vspace*{-4mm}\\
%%%%%%%%%%%%%%%%%%%%%%%%%%%%%%%%%%
$12T57$ & \multirow{5}{*}{
\begin{minipage}{12cm} 
{\rm (I)} $\Sha(T)=0$ if and only if 
there exists a place $v$ of $k$ such that 
${\rm Syl}_2(G)\simeq(C_4\times C_2)\rtimes C_4\leq G_v$;\\
{\rm (II)} $\Sha(T)\simeq\bZ/2\bZ$ if and only if 
there exists a place $v$ of $k$ such that 
$C_4\times C_2\leq G_v$ where 
$|{\rm Orb}_{C_4\times C_2}(i)|=4$ for $4$ $($resp. $8$$)$ 
of $1\leq i\leq 12$ and 
$|{\rm Orb}_{C_4\times C_2}(i)|=2$ for $8$ $($resp. $4$$)$ 
of $1\leq i\leq 12$ 
and {\rm (I)} does not hold
\end{minipage}}\\
$\simeq ((C_4\times C_2)\rtimes C_4)\rtimes C_3$ & \\
 & \\
 & \\
 & \vspace*{1mm}\\\hline
%%%%%%%%%%%%%%%%%%%%%%%%%%%%%%%%%%
\end{tabular}
\end{center}~\\
\vspace*{5mm}

\begin{center}
\vspace*{1mm}
Table $2$-$5$: $\Sha(T)$ for $T=R^{(1)}_{K/k}(\bG_m)$ 
and $G={\rm Gal}(L/k)=12T32$ as in Table $1$\vspace*{2mm}\\
\renewcommand{\arraystretch}{1.05}
\begin{tabular}{ll} 
$G$ & $\Sha(T)\leq (\bZ/2\bZ)^{\oplus 2}$\\\hline\vspace*{-4mm}\\
%%%%%%%%%%%%%%%%%%%%%%%%%%%%%%%%%%
$12T32\simeq (C_2)^4\rtimes C_3$\hspace*{10mm} & \multirow{8}{*}{
\begin{minipage}{12cm}
~\vspace*{-1mm}\\
{\rm (I)} $\Sha(T)=0$ if and only if\\ 
{\rm (I-i)} 
there exists a place $v$ of $k$ such that 
${\rm Syl}_2(G)\simeq (C_2)^4\leq G_v$, or\\
{\rm (I-ii)} 
there exist places $v_1,v_2$ of $k$ such that 
$V_4\simeq V_4^{(1)}\leq G_{v_1}$, $V_4\simeq V_4^{(2)}\leq G_{v_2}$ 
where 
$|{\rm Orb}_{V_4^{(1)}}(i)|=|{\rm Orb}_{V_4^{(2)}}(i)|=2$ 
for any $1\leq i\leq 12$ 
and $V_4^{(1)}\neq V_4^{(2)}$;\\
%and $G_{v_1}\cap G_{v_2}=1$;\\
%$($equivalently, 
%$(\bZ/2\bZ)^{\oplus 2}\leq H^3(G_{v_1},\bZ)\oplus H^3(G_{v_2},\bZ)
%\leq H^3(G,\bZ)\simeq (\bZ/2\bZ)^{\oplus 4}$$)$
{\rm (II)} $\Sha(T)\simeq \bZ/2\bZ$ if and only if 
there exists a place $v$ of $k$ such that 
$V_4\leq G_v$ 
where 
$|{\rm Orb}_{V_4}(i)|=2$ for any $1\leq i\leq 12$ 
and {\rm (I)} does not hold
\end{minipage}}\\
 & \\
 & \\
 & \\
 & \\
 & \\
 & \\
 & \vspace*{1mm}\\\hline
%%%%%%%%%%%%%%%%%%%%%%%%%%%%%%%%%%
\end{tabular}
\end{center}~\\

\newpage
\begin{remark}
As a consequence of Theorem \ref{thmain2}, 
we obtain the 
Tamagawa number $\tau(T)$ of norm one tori $T=R^{(1)}_{k/k}(\bG_m)$ 
over a number field $k$ 
via Ono's formula $\tau(T)=|H^1(G,J_{G/H})|/|\Sha(T)|$ 
where $H^1(G,J_{G/H})$ is given as in 
Section 9 (Appendix) of the arXiv version (arXiv:1910.01469) of \cite{HKY} 
(see Ono \cite[Main theorem, page 68]{Ono63}, 
Voskresenskii \cite[Theorem 2, page 146]{Vos98} and 
Hoshi, Kanai and Yamasaki \cite[Section 8, Application 2]{HKY}). 
%and $\Sha(T)$ is given as in Theorem \ref{thmain2}. 
\end{remark}

%%%%%%%%%%%%%%
We organize this paper as follows. 
In Section \ref{S2}, 
we prepare basic definitions and known results 
about algebraic $k$-tori, in particular, norm one tori.  
In Section \ref{S3}, 
we recall some known results about Hasse norm principle. 
In Section \ref{S4}, we give 
some known results about norm one tori and $k$-rationality. 
In particular, 
we recall a flabby resolution of $G$-lattices 
which is a basic tool to investigate algebraic $k$-tori. 
In Section \ref{S5}, we give the proof of Theorem \ref{thmain1}. 
In Section \ref{S6}, we recall 
Drakokhrust and Platonov's method for 
the Hasse norm principle for $K/k$ 
and results in \cite[Section 6]{HKY}. 
In Section \ref{S7}, we prove Theorem \ref{thmain2} 
by using Drakokhrust and Platonov's method and 
some new useful functions of GAP \cite{GAP}. 
In Section \ref{S8}, 
GAP algorithms are given 
which are also available from 
{\tt https://www.math.kyoto-u.ac.jp/\~{}yamasaki/Algorithm/Norm1ToriHNP}.
\begin{acknowledgments}
We would like to thank 
Ming-chang Kang, Shizuo Endo and Boris Kunyavskii 
for giving us useful and valuable comments. 
We also thank the referee for very careful reading of the manuscript. 
\end{acknowledgments}

%%%%%%%%%%%%%%%%%%%%%%%%%%%%%%%%%%%%%%%%%%%%%%%%%%%%%%%%%%%%%%%%%%%%%%%%%%
\section{Algebraic tori and Norm one tori}\label{S2}
Let $k$ be a field, 
$\overline{k}$ be a fixed separable closure of $k$ and 
$\mathcal{G}={\rm Gal}(\overline{k}/k)$ be the absolute Galois group of $k$. 
Let $T$ be an algebraic $k$-torus, 
i.e. a group $k$-scheme with fiber product (base change) 
$T\times_k \overline{k}=
T\times_{{\rm Spec}\, k}\,{\rm Spec}\, \overline{k}
\simeq (\bG_{m,\overline{k}})^n$; 
$k$-form of the split torus $(\bG_m)^n$. 
Then there exists a finite Galois extension $K/k$ 
with Galois group $G={\rm Gal}(K/k)$ such that 
$T$ splits over $K$: $T\times_k K\simeq (\bG_{m,K})^n$. 
It is also well-known that 
there is the duality between the category of $G$-lattices, 
i.e. finitely generated $\bZ[G]$-modules which are $\bZ$-free 
as abelian groups, 
and the category of algebraic $k$-tori which split over $K$ 
(see Ono \cite[Section 1.2]{Ono61}, 
Voskresenskii \cite[page 27, Example 6]{Vos98} and 
Knus, Merkurjev, Rost and Tignol \cite[page 333, Proposition 20.17]{KMRT98}). 
Indeed, if $T$ is an algebraic $k$-torus, then the character 
module $\widehat{T}={\rm Hom}(T,\bG_m)$ of $T$ 
may be regarded as a $G$-lattice. 
Let $X$ be a smooth $k$-compactification of $T$, 
i.e. smooth projective $k$-variety $X$ 
containing $T$ as a dense open subvariety, 
and $\overline{X}=X\times_k\overline{k}$. 
There exists such a smooth $k$-compactification of an algebraic $k$-torus $T$ 
over any field $k$ (due to Hironaka \cite{Hir64} for ${\rm char}\, k=0$, 
see Colliot-Th\'{e}l\`{e}ne, Harari and Skorobogatov 
\cite[Corollaire 1]{CTHS05} for any field $k$). 

%%%%%%%%%%%%%%%%%%%%%%%%%%%%%%%%%%%
%
Recall the following basic definitions of $G$-lattices. 
%%
%A $\mathcal{G}$-lattice $P$ is said to be {\it permutation} if 
%$P$ has a $\bZ$-basis permuted by $\mathcal{G}$ 
%and 
%a $\mathcal{G}$-lattice $F$ is said to be {\it flabby} 
%(resp. {\it coflabby}) 
%if $\widehat H^{-1}(\mathcal{H},F)=0$ 
%(resp. $H^1(\mathcal{H},F)=0$) 
%for any closed subgroup $\mathcal{H}\leq \mathcal{G}$ 
%where $\widehat H$ is the Tate cohomology. 
\begin{definition}
%[Permutation, stably permutation, invertible, 
%flabby and coflabby $G$-lattices]
Let $G$ be a finite group and $M$ be a $G$-lattice, 
i.e. finitely generated $\bZ[G]$-module which is $\bZ$-free 
as an abelian group.\\
{\rm (i)} $M$ is called a {\it permutation} $G$-lattice 
if $M$ has a $\bZ$-basis permuted by $G$, 
i.e. $M\simeq \oplus_{1\leq i\leq m}\bZ[G/H_i]$ 
for some subgroups $H_1,\ldots,H_m$ of $G$.\\
{\rm (ii)} $M$ is called a {\it stably permutation} $G$-lattice 
if $M\oplus P\simeq P^\prime$ 
for some permutation $G$-lattices $P$ and $P^\prime$.\\
{\rm (iii)} $M$ is called {\it invertible} (or {\it permutation projective}) 
if it is a direct summand of a permutation $G$-lattice, 
i.e. $P\simeq M\oplus M^\prime$ for some permutation $G$-lattice 
$P$ and a $G$-lattice $M^\prime$.\\
{\rm (iv)} $M$ is called {\it flabby} (or {\it flasque}) if $\widehat H^{-1}(H,M)=0$ 
for any subgroup $H$ of $G$ where $\widehat H$ is the Tate cohomology.\\
{\rm (v)} $M$ is called {\it coflabby} (or {\it coflasque}) if $H^1(H,M)=0$
for any subgroup $H$ of $G$.
\end{definition}
\begin{theorem}[{Voskresenskii \cite[Section 4, page 1213]{Vos69}, \cite[Section 3, page 7]{Vos70}, see also \cite[Section 4.6]{Vos98}, \cite[Theorem 1.9]{Kun07}, \cite{Vos74} and \cite[Theorem 5.1, page 19]{CT07} for any field $k$}]\label{thVos69}
Let $k$ be a field 
and $\mathcal{G}={\rm Gal}(\overline{k}/k)$. 
Let $T$ be an algebraic $k$-torus, 
$X$ be a smooth $k$-compactification of $T$ 
and $\overline{X}=X\times_k\overline{k}$. 
Then there exists an exact sequence of $\mathcal{G}$-lattices 
\begin{align*}
0\to \widehat{T}\to \widehat{Q}\to {\rm Pic}\,\overline{X}\to 0%\label{seM}
\end{align*}
where $\widehat{Q}$ is permutation 
and ${\rm Pic}\ \overline{X}$ is flabby. 
\end{theorem}
We have 
$H^1(k,{\rm Pic}\,\overline{X})\simeq H^1(G,{\rm Pic}\, X_K)$ 
where $K$ is the splitting field of $T$, $G={\rm Gal}(K/k)$ and 
$X_K=X\times_k K$. 
Hence Theorem \ref{thVos69} says that 
for $G$-lattices $M=\widehat{T}$ and $P=\widehat{Q}$, 
the exact sequence 
$0\to M\to P\to {\rm Pic}\, X_K\to 0$ 
gives a flabby resolution of 
$M$ and the flabby class of $M$ is 
$[M]^{fl}=[{\rm Pic}\ X_K]$ as $G$-lattices 
(see Theorem \ref{thEM}). 

%%%%%%%%%%%%%%%%%%%%%%%%%%%%%%%%%%%%%%%%%%

\begin{definition}
The kernel of the norm map 
$R_{K/k}(\bG_m)\rightarrow \bG_m$ 
is called {\it the norm one torus 
$R^{(1)}_{K/k}(\bG_m)$ of $K/k$} 
where 
$R_{K/k}$ is the Weil restriction 
(see \cite[page 37, Section 3.12]{Vos98}).
\end{definition}
A norm one torus $R^{(1)}_{K/k}(\bG_m)$ 
is biregularly isomorphic to the norm hyper surface 
$f(x_1,\ldots,x_n)=1$ where 
$f\in k[x_1,\ldots,x_n]$ is the polynomial of total 
degree $n$ defined by the norm map $N_{K/k}:K^\times\to k^\times$.
When $K/k$ is a finite Galois extension, %i.e. $H=1$, 
we have that: 

\begin{theorem}[{Voskresenskii \cite[Theorem 7]{Vos70}, Colliot-Th\'{e}l\`{e}ne and Sansuc \cite[Proposition 1]{CTS77}}]
Let $k$ be a field and 
$K/k$ be a finite Galois extension with Galois group $G={\rm Gal}(K/k)$. 
Let $T=R^{(1)}_{K/k}(\bG_m)$ be the norm one torus of $K/k$ 
and $X$ be a smooth $k$-compactification of $T$. 
Then 
$H^1(H,{\rm Pic}\, X_K)\simeq H^3(H,\bZ)$ for any subgroup $H$ of $G$. 
In particular, 
$H^1(k,{\rm Pic}\, \overline{X})\simeq
H^1(G,{\rm Pic}\, X_K)\simeq H^3(G,\bZ)$ which is isomorphic to 
the Schur multiplier $M(G)$ of $G$.
\end{theorem}
In other words, for $G$-lattice $J_G=\widehat{T}$, 
$H^1(H,[J_G]^{fl})\simeq H^3(H,\bZ)$ for any subgroup $H$ of $G$ 
and $H^1(G,[J_G]^{fl})\simeq H^3(G,\bZ)\simeq H^2(G,\bQ/\bZ)$; 
the Schur multiplier of $G$. 
By the exact sequence $0\to\bZ\to\bZ[G]\to J_G\to 0$, 
we also have $\delta:H^1(G,J_G)\simeq H^2(G,\bZ)\simeq G^{ab}\simeq 
G/[G,G]$ where $\delta$ is the connecting homomorphism and 
$G^{ab}$ is the abelianization of $G$ (cf. Section \ref{S4}). 

%%%%%%%%%%%%%%%%%%%%%%%%%%%%%%%%%%%%%%%%%%%%%
\section{Hasse principle and Hasse norm principle}\label{S3}
%%%%%%%%%%%%%%%%%%%%%%%%%%%%%%%%%%%

%%%%%%%%%%%%%%%%%%%%%%%%%%%%%%%%%%%%%%%%%%%%%%%%%%%%%%%%%%%%5
%Let $k$ be a global field %a number field, 
Let $k$ be a global field, 
i.e. a number field (a finite extension of $\bQ$) 
or a function field of an algebraic curve over 
$\bF_q$ (a finite extension of $\bF_q(t))$.

\begin{definition}
Let $T$ be an algebraic $k$-torus 
and $T(k)$ be the group of $k$-rational points of $T$. 
Then $T(k)$ 
embeds into $\prod_{v\in V_k} T(k_v)$ by the diagonal map 
where %$T(k_v)$ is the group of $k_v$-rational points and 
$V_k$ is the set of all places of $k$ and 
$k_v$ is the completion of $k$ at $v$. 
Let $\overline{T(k)}$ be the closure of $T(k)$  
in the product $\prod_{v\in V_k} T(k_v)$. 
The group 
\begin{align*}
A(T)=\left(\prod_{v\in V_k} T(k_v)\right)/\overline{T(k)}
\end{align*}
is called {\it the kernel of the weak approximation} of $T$. 
We say that {\it $T$ has the weak approximation property} if $A(T)=0$. 
\end{definition}

\begin{definition}
Let $E$ be a principal homogeneous space (= torsor) under $T$.  
{\it Hasse principle holds for $E$} means that 
if $E$ has a $k_v$-rational point for all $k_v$, 
then $E$ has a $k$-rational point. 
The set $H^1(k,T)$ classifies all such torsors $E$ up 
to (non-unique) isomorphism. 
We define {\it the Shafarevich-Tate group} of $T$: 
\begin{align*}
\Sha(T)={\rm Ker}\left\{H^1(k,T)\xrightarrow{\rm res} \bigoplus_{v\in V_k} 
H^1(k_v)\right\}.
\end{align*}
%We say that $T$ {\it satisfies the Hasse principle} 
%if $\Sha(T)=0$, i.e. Hasse principle holds for all torsors $E$ under $T$. 
Then 
Hasse principle holds for all torsors $E$ under $T$ 
if and only if $\Sha(T)=0$. 
\end{definition}

%%%%%%%%%%%%%%%%%%%%%%%%%%%%%%%%%%%%%%%%%%%
%
\begin{theorem}[{Voskresenskii \cite[Theorem 5, page 1213]{Vos69}, 
\cite[Theorem 6, page 9]{Vos70}, see also \cite[Section 11.6, Theorem, page 120]{Vos98}}]\label{thV}
Let $k$ be a global field, 
$T$ be an algebraic $k$-torus and $X$ be a smooth $k$-compactification of $T$. %over $k$. 
Then there exists an exact sequence
\begin{align*}
0\to A(T)\to H^1(k,{\rm Pic}\,\overline{X})^{\vee}\to \Sha(T)\to 0
\end{align*}
where $M^{\vee}={\rm Hom}(M,\bQ/\bZ)$ is the Pontryagin dual of $M$. 
Moreover, if $L$ is the splitting field of $T$ and $L/k$ 
is an unramified extension, then $A(T)=0$ and 
$H^1(k,{\rm Pic}\,\overline{X})^{\vee}\simeq \Sha(T)$. 
\end{theorem}
For the last assertion, see \cite[Section 11.5]{Vos98}. 
It follows that 
$H^1(k,{\rm Pic}\,\overline{X})=0$ if and only if $A(T)=0$ and $\Sha(T)=0$, 
i.e. $T$ has the weak approximation property and 
Hasse principle holds for all torsors $E$ under $T$. 
Theorem \ref{thV} was generalized 
to the case of linear algebraic groups by Sansuc \cite{San81}.

%%%%%%%%%%%%%%%%%%%%%%%%%%%%%%%%%%%%%%%%%%%%%%%%%%%%%%%%%%%%%%%%%
\begin{definition}
Let $k$ be a 
global field, 
$K/k$ be a finite extension and 
$\bA_K^\times$ be the idele group of $K$. 
We say that {\it the Hasse norm principle holds for $K/k$} 
if $(N_{K/k}(\bA_K^\times)\cap k^\times)/N_{K/k}(K^\times)=1$ 
where $N_{K/k}$ is the norm map. 
\end{definition}

Hasse \cite[Satz, page 64]{Has31} proved that 
the Hasse norm principle holds for any cyclic extension $K/k$ 
but does not hold for bicyclic extension $\bQ(\sqrt{-39},\sqrt{-3})/\bQ$. 
For Galois extensions $K/k$, Tate \cite{Tat67} gave the following theorem:
%

%Let $K/k$ be a field extension of degree $n$, 
%$L/k$ be the Galois closure of $K/k$ with 
%$G={\rm Gal}(L/k)=nTm\leq S_n$. 

\begin{theorem}[{Tate \cite[page 198]{Tat67}}]\label{thTate}
Let $k$ be a global field, $K/k$ be a finite Galois extension 
with Galois group ${\rm Gal}(K/k)\simeq G$. 
Let $V_k$ be the set of all places of $k$ 
and $G_v$ be the decomposition group of $G$ at $v\in V_k$. 
Then 
\begin{align*}
(N_{K/k}(\bA_K^\times)\cap k^\times)/N_{K/k}(K^\times)\simeq 
{\rm Coker}\left\{\bigoplus_{v\in V_k}\widehat H^{-3}(G_v,\bZ)\xrightarrow{\rm cores}\widehat H^{-3}(G,\bZ)\right\}
\end{align*}
where $\widehat H$ is the Tate cohomology. 
In particular, the Hasse norm principle holds for $K/k$ 
if and only if the restriction map 
$H^3(G,\bZ)\xrightarrow{\rm res}\bigoplus_{v\in V_k}H^3(G_v,\bZ)$ 
is injective. 
\end{theorem}
If $G\simeq C_n$ is cyclic, then 
$\widehat H^{-3}(G,\bZ)\simeq H^3(G,\bZ)\simeq H^1(G,\bZ)=0$ 
and hence the Hasse's original theorem follows. 
If there exists a place $v$ of $k$ such that $G_v=G$, then 
the Hasse norm principle also holds for $K/k$. 
For example, the Hasse norm principle holds for $K/k$ with 
$G\simeq V_4$ if and only if 
there exists a place $v$ of $k$ such that $G_v=V_4$ because 
$H^3(V_4,\bZ)\simeq\bZ/2\bZ$ and $H^3(C_2,\bZ)=0$. 
The Hasse norm principle holds for $K/k$ with 
$G\simeq (C_2)^3$ if and only if {\rm (i)} 
there exists a place $v$ of $k$ such that $G_v=G$ 
or {\rm (ii)} 
there exist places $v_1,v_2,v_3$ of $k$ such that 
$G_{v_i}\simeq V_4$ and 
$H^3(G,\bZ)\xrightarrow{\rm res}
H^3(G_{v_1},\bZ)\oplus H^3(G_{v_2},\bZ)\oplus H^3(G_{v_3},\bZ)$ 
is isomorphism because 
$H^3(G,\bZ)\simeq(\bZ/2\bZ)^{\oplus 3}$ and 
$H^3(V_4,\bZ)\simeq \bZ/2\bZ$. 

The Hasse norm principle for Galois extensions $K/k$ 
was investigated by Gerth \cite{Ger77}, \cite{Ger78} and 
Gurak \cite{Gur78a}, \cite{Gur78b}, \cite{Gur80} 
(see also \cite[pages 308--309]{PR94}). 
For non-Galois extension $K/k$, 
the Hasse norm principle was investigated by 
Bartels \cite{Bar81a} ($[K:k]=p$; prime), 
\cite{Bar81b} (${\rm Gal}(L/k)\simeq D_n$), 
Voskresenskii and Kunyavskii \cite{VK84} (${\rm Gal}(L/k)\simeq S_n$) and 
Macedo \cite{Mac20} (${\rm Gal}(L/k)\simeq A_n$) 
where $L/k$ be the Galois closure of $K/k$ (see also \cite[Section 1]{HKY}).\\

Ono \cite{Ono63} established the relationship 
between the Hasse norm principle for $K/k$ 
and the Hasse principle for all torsors $E$ under
the norm one torus $R^{(1)}_{K/k}(\bG_m)$. 
%%%%%%%%%%%%%%%%%%%%%%%%%%%%%%%%%%%%%%
The quotient group $k^\times/N_{K/k}(K^\times)$ 
is related to 
the norm one torus $R^{(1)}_{K/k}(\bG_m)$ of $K/k$ 
because by taking the cohomology of the exact sequence 
\begin{align*}
1\to R^{(1)}_{K/k}(\bG_m)\to R_{K/k}(\bG_m)\xrightarrow{N}\bG_m\to 1,
\end{align*}
we have the exact sequence 
\begin{align*}
K^\times\xrightarrow{N_{K/k}}k^\times\to H^1(k,R^{(1)}_{K/k}(\bG_m))\to 
H^1(k,R_{K/k}(\bG_m))\simeq H^1(K,\bG_m)=0\ {\rm (by\ Hilbert\ 90)}
\end{align*} 
and (see Platonov and Rapinchuk \cite[Lemma 2.5, page 73]{PR94})
\begin{align*}
H^1(k,R^{(1)}_{K/k}(\bG_m))\simeq k^\times/N_{K/k}(K^\times). 
\end{align*}
Similarly, the Hasse norm principle for $K/k$ is 
related to the norm one torus $R^{(1)}_{K/k}(\bG_m)$ of $K/k$ 
(see also Platonov and Rapinchuk \cite[Section 6.3]{PR94}): 
\begin{theorem}[{Ono \cite[page 70]{Ono63}, see also Platonov \cite[page 44]{Pla82}, Kunyavskii \cite[Remark 3]{Kun84}, Platonov and Rapinchuk \cite[page 307]{PR94}}]\label{thOno}
Let $k$ be a global field and $K/k$ be a finite extension. 
Then 
\begin{align*}
\Sha(R^{(1)}_{K/k}(\bG_m))\simeq (N_{K/k}(\bA_K^\times)\cap k^\times)/N_{K/k}(K^\times).
\end{align*}
In particular, $\Sha(R^{(1)}_{K/k}(\bG_m))=0$ if and only if 
the Hasse norm principle holds for $K/k$. 
\end{theorem}

\begin{remark}
Applying Theorem \ref{thV} to $T=R^{(1)}_{K/k}(\bG_m)$,  
it follows from Theorem \ref{thOno} that 
$H^1(k,{\rm Pic}\,\overline{X})=0$ if and only if 
$A(T)=0$ and $\Sha(T)=0$, 
i.e. 
$T$ has the weak approximation property and 
the Hasse norm principle holds for $K/k$. 
In the algebraic language, 
the latter condition $\Sha(T)=0$ means that 
for the corresponding norm hyper surface $f(x_1,\ldots,x_n)=b$, 
it has a $k$-rational point 
if and only if it has a $k_v$-rational point 
for any place $v$ of $k$ where 
$f\in k[x_1,\ldots,x_n]$ is the polynomial of total 
degree $n$ defined by the norm map $N_{K/k}:K^\times\to k^\times$ 
and $b\in k^\times$ 
(see \cite[Example 4, page 122]{Vos98}).
\end{remark}

%%%%%%%%%%%%%%%%%%%%%%%%%%%%%%%%%%%
\section{Flabby resolution of $G$-lattice and rationality of algebraic $k$-torus}\label{S4}

Let $k$ be a field, $K/k$ be a separable field extension of degree $n$ 
and $L/k$ be the Galois closure of $K/k$. 
Let $G={\rm Gal}(L/k)$ and $H={\rm Gal}(L/K)$ with $[G:H]=n$. 
The Galois group $G$ may be regarded as a transitive subgroup of 
the symmetric group $S_n$ of degree $n$ via an injection $G\to S_n$ 
which is derived from the action of $G$ on the left cosets 
$\{g_1H,\ldots,g_nH\}$ by $g(g_iH)=(gg_i)H$ for any $g\in G$ 
and we may assume that 
$H$ is the stabilizer of one of the letters in $G$, 
i.e. $L=k(\theta_1,\ldots,\theta_n)$ and $K=k(\theta_i)$ for some 
$1\leq i\leq n$.
%%%%%%%%%%%%%%%%%%%%%%%% 
The norm one torus $R^{(1)}_{K/k}(\bG_m)$ has the 
Chevalley module $J_{G/H}$ as its character module 
and the field $L(J_{G/H})^G$ as its function field 
where $J_{G/H}=(I_{G/H})^\circ={\rm Hom}_\bZ(I_{G/H},\bZ)$ 
is the dual lattice of $I_{G/H}={\rm Ker}\, \varepsilon$ and 
$\varepsilon : \bZ[G/H]\rightarrow \bZ$ is the augmentation map 
(see \cite[Section 4.8]{Vos98}). 
We have the exact sequence 
\begin{align*}
0\rightarrow \bZ\rightarrow \bZ[G/H]\rightarrow J_{G/H}\rightarrow 0
\end{align*}
and rank $J_{G/H}=n-1$. 
Write $J_{G/H}=\oplus_{1\leq i\leq n-1}\bZ u_i$. 
We define the action of $G$ on $L(x_1,\ldots,x_{n-1})$ by 
%\begin{align*}
$\sigma(x_i)=\prod_{j=1}^{n-1} x_j^{a_{i,j}} (1\leq i\leq n-1)$ %\label{acts}
%\end{align*}
for any $\sigma\in G$, when $\sigma (u_i)=\sum_{j=1}^{n-1} a_{i,j} u_j$ 
$(a_{i,j}\in\bZ)$. 
Then the invariant field $L(x_1,\ldots,x_{n-1})^G$ 
may be identified with the function field of the norm 
one torus $R^{(1)}_{K/k}(\bG_m)$ (see \cite[Section 1]{EM75}). 

%%%%
Let $T=R^{(1)}_{K/k}(\bG_m)$ be the norm one torus of $K/k$. 
The rationality problem for norm one tori is investigated 
by \cite{EM75}, \cite{CTS77}, \cite{Hur84}, \cite{CTS87}, 
\cite{LeB95}, \cite{CK00}, \cite{LL00}, \cite{Flo}, \cite{End11}, 
\cite{HY17}, \cite{HY21}, \cite{HHY20} (see also \cite[Section 2]{HKY}). 

We recall some basic facts of flabby (flasque) $G$-lattices
(see \cite{CTS77}, \cite{Swa83}, \cite[Chapter 2]{Vos98}, \cite[Chapter 2]{Lor05}, \cite{Swa10}).

%\begin{lemma}[Lenstra {\cite[Propositions 1.1 and 1.2]{Len74}, see also Swan 
%\cite[Section 8]{Swa83}}]\label{lemSL}
%Let $E$ be an invertible $G$-lattice.\\
%{\rm (i)} $E$ is flabby and coflabby.\\
%{\rm (ii)} If $C$ is a coflabby $G$-lattice, then any short exact sequence
%$0 \rightarrow C \rightarrow N \rightarrow E \rightarrow 0$ splits.
%\end{lemma}
%%%%%%%%%%%%%%%%%%%%%%%%%%

\begin{definition}[{see \cite[Section 1]{EM75}, \cite[Section 4.7]{Vos98}}]
Let $\cC(G)$ be the category of all $G$-lattices. 
Let $\cS(G)$ be the full subcategory of $\cC(G)$ of all permutation $G$-lattices 
and $\cD(G)$ be the full subcategory of $\cC(G)$ of all invertible $G$-lattices.
Let 
\begin{align*}
\cH^i(G)=\{M\in \cC(G)\mid \widehat H^i(H,M)=0\ {\rm for\ any}\ H\leq G\}\ (i=\pm 1)
\end{align*}
be the class of ``$\widehat H^i$-vanish'' $G$-lattices 
where $\widehat H^i$ is the Tate cohomology. 
Then we have the inclusions 
$\cS(G)\subset \cD(G)\subset \cH^i(G)\subset \cC(G)$ $(i=\pm 1)$. 
\end{definition}

\begin{definition}\label{defCM}
We say that two $G$-lattices $M_1$ and $M_2$ are {\it similar} 
if there exist permutation $G$-lattices $P_1$ and $P_2$ such that 
$M_1\oplus P_1\simeq M_2\oplus P_2$. 
We denote the similarity class of $M$ by $[M]$. 
The set of similarity classes $\cC(G)/\cS(G)$ becomes a 
commutative monoid 
(with respect to the sum $[M_1]+[M_2]:=[M_1\oplus M_2]$ 
and the zero $0=[P]$ where $P\in \cS(G)$). 
\end{definition}

%\begin{theorem}[Endo and Miyata {\cite[Theorem 3.3]{EM75}, 
%Endo and Kang {\cite[Theorem 1.4]{EK17}}}]
%Let $G$ be a finite group. 
%Then the following conditions are equivalent:\\
%{\rm (i)} The commutative monoid $\cH^{-1}(G)/\cS(G)$ is a group;\\ 
%%i.e. any element is invertible;\\
%{\rm (ii)} $G=C_n$, $G=D_m$ $(m: odd)$, 
%$G=C_{q^f}\times D_m$ 
%$(q: odd\ prime,\ f\geq 1,\ m: odd,\ {\rm gcd}\{q,m\}=1)$ 
%where $\langle p\rangle=\bF_{q^f}^\times$
%for any prime divisor $p$ of $m$, or 
%$G=Q_{4m}$ $(m: odd)$ where $p\equiv 3\pmod{4}$ for
%any prime divisor $p$ of $m$. %is congruent to $3$ modulo $4$.
%%{\rm (1)} cyclic group, 
%%{\rm (2)} dihedral group of order $2m$ $(m\ {\rm odd})$, 
%%{\rm (3)} a direct product of a cyclic group of order $q^f$, 
%%$q$ an odd prime, $f\geq 1$ 
%%and a dihedral group of order $2m$, $m$ odd, where
%% each prime divisor of $m$ is a primitive $q^{f-1}(q-1)$-th 
%% root of unity modulo $q^f$, or 
%%{\rm (4)} a generalized quaternion group of order $4m$ where 
%%$m$ odd, where each prime divisor of $m$ is congruent to $3$ modulo $4$.
%\end{theorem}
%%
%
\begin{theorem}[{Endo and Miyata \cite[Lemma 1.1]{EM75}, Colliot-Th\'el\`ene and Sansuc \cite[Lemma 3]{CTS77}, 
see also \cite[Lemma 8.5]{Swa83}, \cite[Lemma 2.6.1]{Lor05}}]\label{thEM}
For any $G$-lattice $M$,
there exists a short exact sequence of $G$-lattices
$0 \rightarrow M \rightarrow P \rightarrow F \rightarrow 0$
where $P$ is permutation and $F$ is flabby.
\end{theorem}
\begin{definition}\label{defFlabby}
The exact sequence $0 \rightarrow M \rightarrow P \rightarrow F \rightarrow 0$ 
as in Theorem \ref{thEM} is called a {\it flabby resolution} of the $G$-lattice $M$.
$\rho_G(M)=[F] \in \cC(G)/\cS(G)$ is called {\it the flabby class} of $M$,
denoted by $[M]^{fl}=[F]$.
Note that $[M]^{fl}$ is well-defined: 
if $[M]=[M^\prime]$, $[M]^{fl}=[F]$ and $[M^\prime]^{fl}=[F^\prime]$
then $F \oplus P_1 \simeq F^\prime \oplus P_2$
for some permutation $G$-lattices $P_1$ and $P_2$,
and therefore $[F]=[F^\prime]$ (cf. \cite[Lemma 8.7]{Swa83}). 
We say that $[M]^{fl}$ is {\it invertible} if 
$[M]^{fl}=[E]$ for some invertible $G$-lattice $E$. 
\end{definition}

For $G$-lattice $M$, 
it is not difficult to see 
\begin{align*}
\textrm{permutation}\ \ 
\Rightarrow\ \ 
&\textrm{stably\ permutation}\ \ 
\Rightarrow\ \ 
\textrm{invertible}\ \ 
\Rightarrow\ \ 
\textrm{flabby\ and\ coflabby}\\
&\hspace*{8mm}\Downarrow\hspace*{34mm} \Downarrow\\
&\hspace*{7mm}[M]^{fl}=0\hspace*{10mm}\Rightarrow\hspace*{5mm}[M]^{fl}\ 
\textrm{is\ invertible}.
\end{align*}

The above implications in each step cannot be reversed 
(see, for example, \cite[Section 1]{HY17}). 

\begin{definition}
Let $k$ be a field, and $K$ and $K^\prime$ 
be finitely generated field extensions of $k$.\\
{\rm (i)} $K$ is called {\it rational over $k$} 
(or {\it $k$-rational} for short) 
if $K$ is purely transcendental over $k$, 
i.e. $K$ is isomorphic to $k(x_1,\ldots,x_n)$, 
the rational function field over $k$ with $n$ variables $x_1,\ldots,x_n$ 
for some integer $n$.\\
{\rm (ii)} 
$K$ is called {\it stably $k$-rational} 
if $K(y_1,\ldots,y_m)$ is $k$-rational for some algebraically 
independent elements $y_1,\ldots,y_m$ over $K$.\\
{\rm (iii)}
$K$ and $K^\prime$ are called {\it stably $k$-isomorphic} if 
$K(y_1,\ldots,y_m)\simeq K^\prime(z_1,\ldots,z_n)$ over $k$ 
for some algebraically independent elements $y_1,\ldots,y_m$ over $K$ 
and $z_1,\ldots,z_n$ over $K^\prime$.\\
{\rm (iv)} When $k$ is an infinite field, 
%%%
$K$ is called {\it retract $k$-rational} 
if there is a $k$-algebra $R$ contained in $K$ such that 
$K$ is the quotient field of $R$, and 
the identity map $1_R : R\rightarrow R$ factors through a localized 
polynomial ring over $k$, i.e. there is an element $f\in k[x_1,\ldots,x_n]$, 
which is the polynomial ring over $k$, and there are $k$-algebra 
homomorphisms $\varphi : R\rightarrow k[x_1,\ldots,x_n][1/f]$ 
and $\psi : k[x_1,\ldots,x_n][1/f]\rightarrow R$ satisfying 
$\psi\circ\varphi=1_R$ (cf. \cite{Sal84}).\\
{\rm (v)} $K$ is called {\it $k$-unirational} 
if $k\subset K\subset k(x_1,\ldots,x_n)$ for some integer $n$. 
\end{definition}
It is not difficult to see that 
``$k$-rational'' $\Rightarrow$ ``stably $k$-rational'' $\Rightarrow$ 
``retract $k$-rational'' $\Rightarrow$ ``$k$-unirational''. 

Let $L/k$ be a finite Galois extension with Galois group $G={\rm Gal}(L/k)$ 
and $M$ be a $G$-lattice. 
The flabby class $\rho_G(M)=[M]^{fl}$ 
plays crucial role in the rationality problem for 
$L(M)^G$ as follows (see Voskresenskii's fundamental book \cite[Section 4.6]{Vos98} and Kunyavskii \cite{Kun07}, see also e.g. Swan \cite{Swa83}, 
Kunyavskii \cite[Section 2]{Kun90}, 
Lemire, Popov and Reichstein \cite[Section 2]{LPR06}, 
Kang \cite{Kan12}, Yamasaki \cite{Yam12}, Hoshi and Yamasaki \cite{HY17}, 
Hoshi, Kanai and Yamasaki \cite{HKY}):
\begin{theorem}[{Endo and Miyata \cite{EM73}, Voskresenskii \cite{Vos74}, Saltman \cite{Sal84}}]\label{thEM73}
Let $L/k$ be a finite Galois extension with Galois group $G={\rm Gal}(L/k)$. 
Let $M$ and $M^\prime$ be $G$-lattices. 
Then we have:\\
{\rm (i)} $(${\rm Endo and Miyata} \cite[Theorem 1.6]{EM73}$)$ 
$[M]^{fl}=0$ if and only if $L(M)^G$ is stably $k$-rational;\\
{\rm (ii)} $(${\rm Voskresenskii} \cite[Theorem 2]{Vos74}$)$ 
$[M]^{fl}=[M^\prime]^{fl}$ if and only if $L(M)^G$ and $L(M^\prime)^G$ 
are stably $k$-isomorphic;\\
{\rm (iii)} $(${\rm Saltman} \cite[Theorem 3.14]{Sal84}$)$ 
$[M]^{fl}$ is invertible if and only if $L(M)^G$ is 
retract $k$-rational.
\end{theorem}

It is easy to see that all the $1$-dimensional algebraic $k$-tori $T$, 
i.e. the trivial torus $\bG_m$ and the norm one torus 
$R^{(1)}_{K/k}(\bG_m)$ of $K/k$ with $[K:k]=2$, are $k$-rational. 
Voskresenskii \cite{Vos67} showed that 
all the $2$-dimensional algebraic $k$-tori $T$ are $k$-rational. 
Kunyavskii \cite{Kun90} gave a 
rational (stably rational, retract rational) classification of 
$3$-dimensional algebraic $k$-tori. 
Hoshi and Yamasaki \cite[Theorem 1.9 and Theorem 1.12]{HY17} classified stably/retract rational algebraic $k$-tori of dimension $4$ and $5$ (see also \cite[Section 1]{HKY}).

For norm one tori $T=R^{(1)}_{K/k}(\bG_m)$, 
recall that 
the function field $k(T)$ may be regarded as $L(M)^G$ 
for the character module $M=J_{G/H}$ and hence we have: 
\begin{align*}
[J_{G/H}]^{fl}=0\,
\ \ \Rightarrow\ \ [J_{G/H}]^{fl}\ \textrm{is\ invertible}
\ \ \Rightarrow\ \  H^1(G,[J_{G/H}]^{fl})=0\,
\ \ \Rightarrow\ \  A(T)=0\ \textrm{and}\ \Sha(T)=0
\end{align*}
where the last implication holds over a global field $k$ 
(see also Colliot-Th\'{e}l\`{e}ne and Sansuc \cite[page 29]{CTS77}). 
The last conditions mean that %are equivalent to the condition that 
$T$ has the weak approximation property and 
the Hasse norm principle holds for $K/k$ (see Section \ref{S1}). 
In particular, it follows that 
$[J_{G/H}]^{fl}$ is invertible and hence 
$A(T)=0$ and $\Sha(T)=0$ when $G=pTm\leq S_p$ is a transitive 
subgroup of $S_p$ of prime degree $p$ 
and $H=G\cap S_{p-1}\leq G$ with $[G:H]=p$ (see 
Colliot-Th\'{e}l\`{e}ne and Sansuc \cite[Proposition 9.1]{CTS87} and 
\cite[Lemma 2.17]{HY17}). 
Hence the Hasse norm principle holds for $K/k$ when $[K:k]=p$ 
(see the second paragraph after Theorem \ref{thTate} 
and also the first paragraph of Section \ref{S5}).  

A necessary and sufficient condition for the classification 
of stably/retract rational norm one tori $T=R^{(1)}_{K/k}(\bG_m)$ 
with $[K:k]=n\leq 15$, 
but with one exception 
$G=9T27\simeq \PSL_2(\bF_8)$ for the stable rationality, 
was given in 
Hoshi and Yamasaki \cite{HY21} (for the case $n$ is a prime number or 
the case $n\leq 10$)
and Hasegawa, Hoshi and Yamasaki \cite{HHY20} (for $n=12,14,15$). 

%%%%%%%%%%%%%%%%%%%%%%%%%%%%%%%%%
\section{{Proof of Theorem \ref{thmain1}}}\label{S5}

Let $K/k$ be a separable field extension of degree $12$ 
and $L/k$ be the Galois closure of $K/k$. 
Let $G={\rm Gal}(L/k)=12Tm$ $(1\leq m\leq 301)$ 
be a transitive subgroup of $S_{12}$ 
and $H={\rm Gal}(L/K)$ with $[G:H]=12$. 
Let $T=R^{(1)}_{K/k}(\bG_m)$ be the norm one torus of $K/k$ 
of dimension $11$. 
We have the character module $\widehat{T}=J_{G/H}$ of $T$ 
and then 
$H^1(k,{\rm Pic}\,\overline{X})\simeq H^1(G,[J_{G/H}]^{fl})$ 
(see Section \ref{S4}). 
We may assume that 
$H$ is the stabilizer of one of the letters in $G$, 
i.e. $L=k(\theta_1,\ldots,\theta_{12})$ and $K=k(\theta_i)$ for some 
$1\leq i\leq 12$. 
In order to compute $H^1(G,[J_{G/H}]^{fl})$, 
we apply the functions 
{\tt Norm1TorusJ($12,m$)} and 
{\tt FlabbyResolution($G$).actionF} 
for $1\leq m\leq 299$ as in \cite[Chapter 5]{HY17}.  
The function 
{\tt Norm1TorusJ($12,m$)} returns 
$J_{G/H}$ for $G=12Tm\leq S_{12}$ 
and $H$ is the stabilizer of one of the letters in $G$, 
and {\tt FlabbyResolution($G$).actionF} returns a 
flabby class $F=[J_{G/H}]^{fl}$ for $G=12Tm\leq S_{12}$. 
The assertion for $G=12Tm$ $(1\leq m\leq 299)$ 
follows from these computations as in Example \ref{exH1F}. 

For $G=12T300\simeq A_{12}$, $12T301\simeq S_{12}$, 
we did not get an answer by the computer calculations 
because it needs much time and computer resources (memory) 
in computations. 
However, we already know that 
$H^1(G,[J_{G/H}]^{fl})=0$ for 
$G=12T300\simeq A_{12}$, $12T301\simeq S_{12}$ 
by Macedo \cite{Mac20} 
and Voskresenskii and Kunyavskii \cite{VK84} respectively 
(see also \cite[Theorem 4, Corollary]{Vos88} 
and \cite[Theorem 1.10 and Theorem 1.11]{HKY}) . 
%{\it Proof of Theorem \ref{thmain1}}. 

The last assertion follows from Theorem \ref{thV}.\qed\\

Some related functions for Example \ref{exH1F}  
are available from\\ 
%\begin{center}
{\tt https://www.math.kyoto-u.ac.jp/\~{}yamasaki/Algorithm/RatProbNorm1Tori/}.
%\end{center}

\smallskip
\begin{example}[{Computation of $H^1(G,[J_{G/H}]^{fl})$ with $G=12Tm$ $(1\leq m\leq 299)$}]\label{exH1F}
~{}\vspace*{-4mm}\\
{\small 

}
\end{example}

%
%%%%%%%%%%%%%%%%%%%%%%%%%%%%%%%%%
\section{Drakokhrust and Platonov's method}\label{S6}

Let $k$ be a number field, $K/k$ be a finite extension, 
$\bA_K^\times$ be the idele group of $K$ and 
$L/k$ be the Galois closure of $K/k$. 
Let $G={\rm Gal}(L/k)=nTm$ be a transitive subgroup of $S_n$ 
and $H={\rm Gal}(L/K)$ with $[G:H]=n$. 

For $x,y\in G$, we denote $[x,y]=x^{-1}y^{-1}xy$ the commutator of 
$x$ and $y$, and $[G,G]$ the commutator group of $G$. 
Let $V_k$ be the set of all places of $k$ 
and $G_v$ be the decomposition group of $G$ at $v\in V_k$. 

\begin{definition}[{Drakokhrust and Platonov \cite[page 350]{PD85a}, \cite[page 300]{DP87}}]
Let $k$ be a number field, 
$L\supset K\supset k$ be a tower of finite extensions 
where $L$ is normal over $k$. 

We call the group 
\begin{align*}
{\rm Obs}(K/k)=(N_{K/k}(\bA_K^\times)\cap k^\times)/N_{K/k}(K^\times)
\end{align*}
{\it the total obstruction to the Hasse norm principle for $K/k$} 
and 
\begin{align*}
{\rm Obs}_1(L/K/k)=\left(N_{K/k}(\bA_K^\times)\cap k^\times\right)/\left((N_{L/k}(\bA_L^\times)\cap k^\times)N_{K/k}(K^\times)\right)
\end{align*}
{\it the first obstruction to the Hasse norm principle for $K/k$ 
%with respect to 
corresponding to the tower 
$L\supset K\supset k$}. 
\end{definition}

Note that (i) 
${\rm Obs}(K/k)=1$ if and only if 
the Hasse norm principle holds for $K/k$; 
and (ii) ${\rm Obs}_1(L/K/k)
={\rm Obs}(K/k)/(N_{L/k}(\bA_L^\times)\cap k^\times)$. 

Drakokhrust and Platonov gave a formula 
for computing the first obstruction ${\rm Obs}_1(L/K/k)$: 
%for $K/k$: 

\begin{theorem}[{Drakokhrust and Platonov \cite[page 350]{PD85a}, \cite[pages 789--790]{PD85b}, \cite[Theorem 1]{DP87}}]\label{thDP2}
%Let $k$ be a number field, $K/k$ be a finite extension 
%and $L/k$ be the Galois closure of $K/k$.
Let $k$ be a number field, 
$L\supset K\supset k$ be a tower of finite extensions 
where 
$L$ is normal over $k$.  
Let $G={\rm Gal}(L/k)$ and $H={\rm Gal}(L/K)$. %with $[G:H]=n$. 
Then 
\begin{align*}
%\Sha(R^{(1)}_{K/k}(\bG_m))
{\rm Obs}_1(L/K/k)\simeq 
{\rm Ker}\, \psi_1/\varphi_1({\rm Ker}\, \psi_2)
\end{align*}
where 
\begin{align*}
\begin{CD}
H/[H,H] @>\psi_1 >> G/[G,G]\\
@AA\varphi_1 A @AA\varphi_2 A\\
\displaystyle{\bigoplus_{v\in V_k}\left(\bigoplus_{w\mid v} H_w/[H_w,H_w]\right)} @>\psi_2 >> 
\displaystyle{\bigoplus_{v\in V_k} G_v/[G_v,G_v]}, 
\end{CD}
\end{align*}
$\psi_1$, $\varphi_1$ and $\varphi_2$ are defined 
by the inclusions $H\subset G$, $H_w\subset H$ and $G_v\subset G$ respectively, and 
\begin{align*}
\psi_2(h[H_{w},H_{w}])=x^{-1}hx[G_v,G_v]
\end{align*}
for $h\in H_{w}=H\cap x^{-1}hx[G_v,G_v]$ $(x\in G)$.
\end{theorem}

%%%%%%%%%%%%%%%%%%%%%%%%%%%%%

Let $\psi_2^{v}$ be the restriction of $\psi_2$ to the subgroup 
$\bigoplus_{w\mid v} H_w/[H_w,H_w]$ with respect to $v\in V_k$ 
and $\psi_2^{\rm nr}$ (resp. $\psi_2^{\rm r}$) be 
the restriction of $\psi_2$ to the unramified (resp. the ramified) 
places $v$ of $k$. 
\begin{proposition}[{Drakokhrust and Platonov \cite{DP87}}]\label{propDP}
Let $k$, 
$L\supset K\supset k$, 
$G$ and $H$ be as in Theorem \ref{thDP2}.\\
{\rm (i)} $($\cite[Lemma 1]{DP87}$)$ 
Places $w_i\mid v$ of $K$ are in one-to-one correspondence 
with the set of double cosets in the decomposition 
$G=\cup_{i=1}^{r_v} Hx_iG_v$ where $H_{w_i}=H\cap x_iG_vx_i^{-1}$;\\
{\rm (ii)} $($\cite[Lemma 2]{DP87}$)$ 
If $G_{v_1}\leq G_{v_2}$, then $\varphi_1({\rm Ker}\,\psi_2^{v_1})\subset \varphi_1({\rm Ker}\,\psi_2^{v_2})$;\\
{\rm (iii)} $($\cite[Theorem 2]{DP87}$)$ 
$\varphi_1({\rm Ker}\,\psi_2^{\rm nr})=\Phi^G(H)/[H,H]$ 
where $\Phi^G(H)=\langle [h,x]\mid h\in H\cap xHx^{-1}, x\in G\rangle$;\\
{\rm (iv)} $($\cite[Lemma 8]{DP87}$)$ If $[K:k]=p^r$ $(r\geq 1)$ 
and ${\rm Obs}(K_p/k_p)=1$ where $k_p=L^{G_p}$, $K_p=L^{H_p}$, 
$G_p$ and $H_p\leq H\cap G_p$ are $p$-Sylow subgroups of $G$ and $H$ 
respectively, then ${\rm Obs}(K/k)=1$.
\end{proposition}

Note that the inverse direction of Proposition \ref{propDP} (iv) 
does not hold in general. 
For example, if $n=8$, $G=8T13\simeq A_4\times C_2$ and 
there exists a place $v$ of $k$ such that $G_v\simeq V_4$, 
then ${\rm Obs}(K/k)=1$ but $G_2=8T3\simeq (C_2)^3$ 
and ${\rm Obs}(K_2/k_2)\neq 1$ may occur 
(see \cite[Theorem 1.18]{HKY}).

\begin{theorem}[{Drakokhrust and Platonov \cite[Theorem 3, Corollary 1]{DP87}}]\label{thDP87}
Let $k$, 
$L\supset K\supset k$, 
$G$ and $H$ be as in Theorem \ref{thDP2}. 
Let $H_i\leq G_i\leq G$ $(1\leq i\leq m)$, 
$H_i\leq H\cap G_i$, 
$k_i=L^{G_i}$ and $K_i=L^{H_i}$. 
If ${\rm Obs}(K_i/k_i)=1$ for all $1\leq i\leq m$ and 
\begin{align*}
\bigoplus_{i=1}^m \widehat{H}^{-3}(G_i,\bZ)\xrightarrow{\rm cores} 
\widehat{H}^{-3}(G,\bZ)
\end{align*}
is surjective, 
then ${\rm Obs}(K/k)={\rm Obs}_1(L/K/k)$. 
In particular, 
if $[K:k]=n$ is square-free, 
then ${\rm Obs}(K/k)={\rm Obs}_1(L/K/k)$.
\end{theorem}

We note that if $L/k$ is an unramified extension, 
then $A(T)=0$ and $H^1(G,[J_{G/H}]^{fl})\simeq \Sha(T)\simeq 
{\rm Obs}(K/k)$ where $T=R^{(1)}_{K/k}(\bG_m)$ 
(see Theorem \ref{thV} and Theorem \ref{thOno}). 
If, in addition, ${\rm Obs}(K/k)={\rm Obs}_1(L/K/k)$ 
(e.g. $[K:k]=6,10,14,15$; square-free, see Theorem \ref{thDP87}), 
then ${\rm Obs}(K/k)={\rm Obs}_1(L/K/k)=
{\rm Ker}\, \psi_1/\varphi_1({\rm Ker}\, \psi_2^{\rm nr})\simeq 
{\rm Ker}\, \psi_1/(\Phi^G(H)/[H,H])$ 
(see Proposition \ref{propDP} (iii)). 

\begin{theorem}[{Drakokhrust \cite[Theorem 1]{Dra89}, see also Opolka \cite[Satz 3]{Opo80}}]\label{thDra89}
Let $k$, 
$L\supset K\supset k$, 
$G$ and $H$ be as in Theorem \ref{thDP2}. 
Assume that $\widetilde{L}\supset L\supset k$ is 
a tower of Galois extensions with 
$\widetilde{G}={\rm Gal}(\widetilde{L}/k)$ 
and $\widetilde{H}={\rm Gal}(\widetilde{L}/K)$ 
which correspond to a central extension 
$1\to A\to \widetilde{G}\to G\to 1$ with 
$A\cap[\widetilde{G},\widetilde{G}]\simeq M(G)=H^2(G,\bC^\times)$; 
the Schur multiplier of $G$ 
$($this is equivalent to 
the inflation 
$M(G)\to M(\widetilde{G})$ being the zero map, 
see {\rm Beyl and Tappe \cite[Proposition 2.13, page 85]{BT82}}$)$. 
Then 
${\rm Obs}(K/k)={\rm Obs}_1(\widetilde{L}/K/k)$. 
In particular, if $\widetilde{G}$ is a Schur cover of $G$, 
i.e. $A\simeq M(G)$, then ${\rm Obs}(K/k)={\rm Obs}_1(\widetilde{L}/K/k)$. 
\end{theorem}

Indeed, Drakokhrust \cite[Theorem 1]{Dra89} shows that 
${\rm Obs}(K/k)\simeq 
{\rm Ker}\, \widetilde{\psi}_1/\widetilde{\varphi}_1({\rm Ker}\, \widetilde{\psi}_2)$ where the maps $\widetilde{\psi}_1, \widetilde{\psi}_2$ and $\widetilde{\varphi}_1$ are defined as in 
\cite[page 31, the paragraph before Proposition 1]{Dra89}. 
The proof of \cite[Proposition 1]{Dra89} shows that 
this group is the same as ${\rm Obs}_1(\widetilde{L}/K/k)$ 
(see also \cite[Lemma 2, Lemma 3 and Lemma 4]{Dra89}).\\

%%%%%%%%%%%%%%%%%%%%%%%%%%%%%%%%%%%%%%%%%%%%%%%%%%%%%%%%%%%%%%%%%%%%%%
Hoshi, Kanai and Yamasaki \cite[Section 6]{HKY} 
made the following functions of GAP (\cite{GAP}) which 
will be used in the proof of Theorem \ref{thmain2}. \\

{\tt FirstObstructionN($G,H$).ker} returns 
the list $[l_1, [l_2, l_3]]$ where 
$l_1$ is the abelian invariant of the numerator of the first obstruction 
${\rm Ker}\, \psi_1=\langle y_1,\ldots,y_t\rangle$ 
with respect to $G$, $H$ as in Theorem \ref{thDP2}, 
$l_2=[e_1,\ldots,e_m]$ 
is the abelian invariant of 
$H^{ab}=H/[H,H]=\langle x_1,\ldots,x_m\rangle$ with $e_i={\rm order}(x_i)$ 
and 
$l_3=[l_{3,1},\ldots,l_{3,t}]$, 
$l_{3,i}=[r_{i,1},\ldots,r_{i,m}]$ is the list with 
$y_i=x_1^{r_{i,1}}\cdots x_m^{r_{i,m}}$ 
for $H\leq G\leq S_n$.

{\tt FirstObstructionN($G$).ker} returns the same 
as {\tt FirstObstructionN($G,H$).ker} where $H={\rm Stab}_1(G)$ 
is the stabilizer of $1$ in $G\leq S_n$.\\

{\tt FirstObstructionDnr($G,H$).Dnr} returns 
the list $[l_1, [l_2, l_3]]$ where 
$l_1$ is the abelian invariant of the unramified part 
of the denominator of the first obstruction 
$\varphi_1({\rm Ker}\, \psi_2^{\rm nr})=\Phi^G(H)/[H,H]=\langle y_1,\ldots,y_t\rangle$ 
with respect to $G$, $H$ as in Proposition \ref{propDP} (iii), 
$l_2=[e_1,\ldots,e_m]$ 
is the abelian invariant of 
$H^{ab}=H/[H,H]=\langle x_1,\ldots,x_m\rangle$ with $e_i={\rm order}(x_i)$ 
and 
$l_3=[l_{3,1},\ldots,l_{3,t}]$, 
$l_{3,i}=[r_{i,1},\ldots,r_{i,m}]$ is the list with 
$y_i=x_1^{r_{i,1}}\cdots x_m^{r_{i,m}}$ 
for $H\leq G\leq S_n$.

{\tt FirstObstructionDnr($G$).Dnr} returns the same 
as {\tt FirstObstructionDnr($G,H$).Dnr} where $H={\rm Stab}_1(G)$ 
is the stabilizer of $1$ in $G\leq S_n$.\\

{\tt FirstObstructionDr($G,G_v,H$).Dr} returns 
the list $[l_1, [l_2, l_3]]$ where 
$l_1$ is the abelian invariant of the ramified part 
of the denominator of the first obstruction 
$\varphi_1({\rm Ker}\, \psi_2^v)=\langle y_1,\ldots,y_t\rangle$ 
with respect to $G$, $G_v$, $H$ as in Theorem \ref{thDP2}, 
$l_2=[e_1,\ldots,e_m]$ 
is the abelian invariant of 
$H^{ab}=H/[H,H]=\langle x_1,\ldots,x_m\rangle$ with $e_i={\rm order}(x_i)$ 
and 
$l_3=[l_{3,1},\ldots,l_{3,t}]$, 
$l_{3,i}=[r_{i,1},\ldots,r_{i,m}]$ is the list with 
$y_i=x_1^{r_{i,1}}\cdots x_m^{r_{i,m}}$ 
for $G_v, H\leq G\leq S_n$.

{\tt FirstObstructionDr($G,G_v$).Dr} returns the same 
as {\tt FirstObstructionDr($G,G_v,H$).Dr} where $H={\rm Stab}_1(G)$ 
is the stabilizer of $1$ in $G\leq S_n$.\\

{\tt SchurCoverG($G$).SchurCover} 
(resp. {\tt SchurCoverG($G$).epi}) 
returns 
one of the Schur covers $\widetilde{G}$ of $G$ 
(resp. the surjective map $\pi$) 
in a central extension 
$1\to A\to \widetilde{G}\xrightarrow{\pi} G\to 1$ 
with $A\simeq M(G)$; Schur multiplier of $G$ 
(see Karpilovsky \cite[page 16]{Kap87}). 
The Schur covers $\widetilde{G}$ are 
stem extensions, i.e. $A\leq Z(\widetilde{G})\cap [\widetilde{G},\widetilde{G}]$, of the maximal size. 
This function is based on the built-in function 
{\tt EpimorphismSchurCover} in GAP.\\

{\tt MinimalStemExtensions($G$)[$j$].MinimalStemExtension} 
(resp. {\tt MinimalStemExtensions($G$)[$j$].epi}) 
returns 
the $j$-th minimal stem extension $\overline{G}=\widetilde{G}/A^\prime$, 
i.e. $\overline{A}\leq Z(\overline{G})\cap [\overline{G},\overline{G}]$, 
of $G$ provided by the Schur cover $\widetilde{G}$ of $G$ 
via {\tt SchurCoverG($G$).SchurCover} 
where $A^\prime$ is the $j$-th maximal subgroup of $A=M(G)$
(resp. the surjective map $\overline{\pi}$) 
in the commutative diagram 
\begin{align*}
\begin{CD}
1 @>>> A=M(G) @>>> \widetilde{G}@>\pi>> G@>>> 1\\
  @. @VVV @VVV @|\\
1 @>>> \overline{A}=A/A^\prime @>>> \overline{G}=\widetilde{G}/A^\prime @>\overline{\pi}>> G@>>> 1
\end{CD}
\end{align*}
(see Robinson \cite[Exercises 11.4]{Rob96}). 
This function is based on the built-in function 
{\tt EpimorphismSchurCover} in GAP.\\
%%%%%%%%%%%%%%%%%%%%%%%%%%

{\tt ResolutionNormalSeries(LowerCentralSeries($G$),} {\tt$n+1$)} 
(resp. {\tt ResolutionNormalSeries(DerivedSeries} {\tt ($G$),$n+1$)}, 
{\tt ResolutionFiniteGroup($G$,$n+1$)}) returns 
a free resolution $RG$ of $G$ 
when $G$ is nilpotent (resp. solvable, finite). 
This function is the built-in function of 
HAP (\cite{HAP}) in GAP (\cite{GAP}) .\\

%%%%%%%%%%%%%%%%%%%%%%%%%%%%%%%%%
{\tt ResHnZ($RG,RH,n$).HnGZ} (resp. {\tt ResHnZ($RG,RH,n$).HnHZ}) returns 
the abelian invariants of $H^n(G,\bZ)$ (resp. $H^n(H,\bZ)$) 
with respect to Smith normal form, 
for free resolutions $RG$ and $RH$ of $G$ and $H$ respectively.\\

{\tt ResHnZ($RG,RH,n$).Res} returns 
the list $L=[l_1,\ldots,l_s]$ where 
$H^n(G,\bZ)=\langle x_1,\ldots,x_s\rangle\xrightarrow{\rm res} 
H^n(H,\bZ)=\langle y_1,\ldots,y_t\rangle$, 
${\rm res}(x_i)=\prod_{j=1}^t y_j^{l_{i,j}}$ 
and $l_i=[l_{i,1},\ldots,l_{i,t}]$ 
for free resolutions $RG$ and $RH$ of $G$ and $H$ respectively.\\

{\tt ResHnZ($RG,RH,n$).Ker} returns 
the list $L=[l_1,[l_2,l_3]]$ 
where $l_1$ is the abelian invariant of 
${\rm Ker}\{H^n(G,\bZ)$ $\xrightarrow{\rm res}$ 
$H^n(H,\bZ)\}=\langle y_1,\ldots,y_t\rangle$, 
$l_2=[d_1,\ldots,d_s]$ is the abelian invariant of $H^n(G,\bZ)=\langle x_1,\ldots,x_s\rangle$ with $d_i={\rm ord}(x_i)$
and $l_3=[l_{3,1},\ldots,l_{3,t}]$, 
$l_{3,j}=[r_{j,1},\ldots,r_{j,s}]$ is the list with 
$y_j=x_1^{r_{j,1}}\cdots x_s^{r_{j,s}}$ 
for free resolutions $RG$ and $RH$ of $G$ and $H$ respectively.\\

{\tt ResHnZ($RG,RH,n$).Coker} returns 
the list $L=[l_1,[l_2,l_3]]$ 
where $l_1=[e_1,\ldots,e_t]$ is the abelian invariant of 
${\rm Coker}\{H^n(G,\bZ)$ $\xrightarrow{\rm res}$ 
$H^n(H,\bZ)\}=\langle \overline{y_1},\ldots,\overline{y_t}\rangle$ 
with $e_j={\rm ord}(\overline{y_j})$, 
$l_2=[d_1,\ldots,d_s]$ is the abelian invariant of $H^n(H,\bZ)=\langle x_1,\ldots,x_s\rangle$ with $d_i={\rm ord}(x_i)$ 
and $l_3=[l_{3,1},\ldots,l_{3,t}]$, 
$l_{3,j}=[r_{j,1},\ldots,r_{j,s}]$ is the list with 
$\overline{y_j}=\overline{x_1}^{r_{j,1}}\cdots \overline{x_s}^{r_{j,s}}$ 
for free resolutions $RG$ and $RH$ of $G$ and $H$ respectively.\\

{\tt KerResH3Z($G,H$)} returns the list $L=[l_1,[l_2,l_3]]$ 
where $l_1$ is the abelian invariant of 
${\rm Ker}\{H^3(G,\bZ)\xrightarrow{\rm res}\oplus_{i=1}^{m^\prime} H^3(G_i,\bZ)\}=\langle y_1,\ldots,y_t\rangle$ 
where $H_i\leq G_i\leq G$, $H_i\leq H\cap G_i$, $[G_i:H_i]=n$ 
and 
the action of $G_i$ on $\bZ[G_i/H_i]$ may be regarded as $nTm$ 
$(n\leq 15, n\neq 12)$ which is not in \cite[Table $1$]{HKY}, 
$l_2=[d_1,\ldots,d_s]$ is the abelian invariant of $H^3(G,\bZ)=\langle x_1,\ldots,x_s\rangle$ with $d_{i^\prime}={\rm ord}(x_{i^\prime})$ 
and $l_3=[l_{3,1},\ldots,l_{3,t}]$, 
$l_{3,j}=[r_{j,1},\ldots,r_{j,s}]$ is the list with 
$y_j=x_1^{r_{j,1}}\cdots x_s^{r_{j,s}}$ 
for groups $G$ and $H$ (cf. Theorem \ref{thDra89}). \\ 

The functions above are available from\\
{\tt https://www.math.kyoto-u.ac.jp/\~{}yamasaki/Algorithm/Norm1ToriHNP}.\\

%%%%%%%%%%%%%%%%
\smallskip
%
%%%%%%%%%%%%%%%%%%%%%%%%%%%%%%%%%
\section{{Proof of Theorem \ref{thmain2}}}\label{S7}

In order to prove Theorem \ref{thmain2}, 
we made the following additional GAP (\cite{GAP}) functions:\\

{\tt KerResH3Z($G,H$:iterator)} returns the same as 
{\tt KerResH3Z($G,H$)} (in Section \ref{S6}) but using 
the built-in function {\tt IteratorByFunctions} of GAP 
in order to run fast (by applying the new 
function {\tt ChooseGiIterator}$(G,H)$ to choose suitable $G_i\leq G$ 
instead of the old one {\tt ChooseGi}$(G,H)$) and also using 
$nTm$ $(n\leq 15)$ which is not in \cite[Table $1$]{HKY} or Table $1$
instead of 
$nTm$ 
$(n\leq 15, n\neq 12)$ which is not in \cite[Table $1$]{HKY} 
(cf. Theorem \ref{thDra89}). \\

{\tt ConjugacyClassesSubgroupsNGHOrbitRep(ConjugacyClassesSubgroups($G$),$H$)} 
returns the list $L=[l_1,\ldots,l_t]$ 
where $t$ is the number of subgroups of $G$ up to conjugacy, 
$l_r=[l_{r,1},\ldots,l_{r,u_r}]$ $(1\leq r\leq t)$, 
$l_{r,s}$ $(1\leq s\leq u_r)$
is a representative of the orbit 
${\rm Orb}_{N_G(H)\backslash G/N_G(G_{v_{r,s}})}(G_{v_{r,s}})$ 
of $G_{v_{r,s}}\leq G$ under the conjugate action of $G$ 
which corresponds to the double coset 
$N_G(H)\backslash G/N_G(G_{v_{r,s}})$ 
with ${\rm Orb}_{G/N_G(G_{v_r})}(G_{v_r})
=\bigcup_{s=1}^{u_r}{\rm Orb}_{N_G(H)\backslash G/N_G(G_{v_{r,s}})}(G_{v_{r,s}})$ 
corresponding to $r$-th subgroup $G_{v_r}\leq G$ up to conjugacy.\\

{\tt MinConjugacyClassesSubgroups($l$)} 
returns the minimal elements of the list $l=\{H_i^G\}$ 
where $H_i^G=\{x^{-1}H_ix\mid x\in G\}$ %$(H_i\leq G)$ 
for some subgroups $H_i\leq G$ which 
satisfy the condition that if $H_i^G, H_r^G\in l$ and $H_i\leq H_r$, 
then $H_j^G\in l$ for any $H_i\leq H_j\leq H_r$.\\

{\tt IsInvariantUnderAutG($l$)} returns {\tt true} if 
the list $l=\{H_i^G\}$ is closed under the action of all 
the automorphisms ${\rm Aut}(G)$ of $G$ 
where $H_i^G=\{x^{-1}H_ix\mid x\in G\}$ $(H_i\leq G)$. 
If not, this returns {\tt false}.\\

The functions above are also available from\\
{\tt https://www.math.kyoto-u.ac.jp/\~{}yamasaki/Algorithm/Norm1ToriHNP}.\\

%%%%%%%%%%%%%%%%%%%%%%%%%%%%%%%%%%%%%%%%%%%%%%%%%%%%%%%%%%%%%%
{\it Proof of Theorem \ref{thmain2}.} 

Let $G={\rm Gal}(L/k)=12Tm\leq S_{12}$ be 
the $m$-th transitive subgroup of $S_{12}$ and 
$H={\rm Gal}(L/K)\leq G$ with $[G:H]=12$. 
Let $V_k$ be the set of all places of $k$ 
and $G_v$ be the decomposition group of $G$ at $v\in V_k$. 

We may assume that 
$H={\rm Stab}_1(G)$ is the stabilizer of $1$ in $G$, 
i.e. $L=k(\theta_1,\ldots,\theta_{12})$ and $K=L^H=k(\theta_1)$. 
Note that  (the multi-set) $\{{\rm Orb}_{G^\prime}(i)\mid 1\leq i\leq 12\}$ 
$(G^\prime\leq G)$ is invariant under the conjugacy actions of $G$, 
i.e. inner automorphisms of $G$.

%By Theorem \ref{thmain1} and Theorem \ref{thV}, 
%it is enough to give a necessary and sufficient condition for $\Sha(T)=0$. 
We also include the proof for Galois cases $L=K$, i.e. $H=1$, 
(see Tate \cite{Tat67} and Theorem \ref{thTate}). 

%%%%%%%%%%%%%%%%%%%%%%%%%%%%%%%%%%%%%%%%%%%%%%%%%%%%%

{\rm (1)} Table $2$-$1$: $G=12Tm$ 
$(m=2,3,4,7,9,16,18,20,33,34,40,43,47,52,55,64,65,70,71,74,75,96,97,
117,122$, $130,132,133,168,171,172,174,176,179,194,232,234,246,261,280)$ 
$($$40$ {\rm cases}$)$ with $\Sha(T)\leq \bZ/2\bZ$ 
$($cf. Tate \cite[page 198]{Tat67} for Galois cases 
$12T2\simeq C_6\times C_2$, 
$12T3\simeq D_6$ and 
$12T4\simeq A_4(12)$$)$. 

We split the $40$ cases into $3$ parts (1-1), (1-2), (1-3) 
according to the method to prove the assertion: 

{\rm (1-1)} 
$m=7,9,34,47,52,55,64,65,74,75,96,97,122,172,174,232,246$ ($17$ cases). 

Applying 
{\tt FirstObstructionN($G$)} and {\tt FirstObstructionDnr($G$)}, 
we obtain that 
${\rm Ker}\, \psi_1/\varphi_1({\rm Ker}\, \psi_2^{\rm nr})$ 
$\simeq$ $\bZ/2\bZ$. 
This implies that ${\rm Obs}(K/k)={\rm Obs}_1(L/K/k)$ 
if $L/k$ is unramified. 
Use 
Theorem \ref{thDP87}. 
Applying 
{\tt KerResH3Z($G,H$:iterator)}, 
we see ${\rm Ker}\{H^3(G,\bZ)\xrightarrow{\rm res}\oplus_{i=1}^{m^\prime} H^3(G_i,\bZ)\}=0$ and hence 
$\oplus_{i=1}^{m^\prime} \widehat{H}^{-3}(G_i,\bZ)\xrightarrow{\rm cores} 
\widehat{H}^{-3}(G,\bZ)$ is surjective.  
It follows from Theorem \ref{thDP87} that 
${\rm Obs}(K/k)={\rm Obs}_1(L/K/k)$. 
Apply the function {\tt FirstObstructionDr($G,G_{v_{r,s}}$)} 
to representatives of 
%%%%%%%%%%%%%%%%%
the orbit 
${\rm Orb}_{N_G(H)\backslash G/N_G(G_{v_{r,s}})}(G_{v_{r,s}})$ 
of $G_{v_{r,s}}\leq G$ under the conjugate action of $G$ 
which corresponds to the double coset 
$N_G(H)\backslash G/N_G(G_{v_{r,s}})$ 
with ${\rm Orb}_{G/N_G(G_{v_r})}(G_{v_r})$
$=$$\bigcup_{s=1}^{u_r}{\rm Orb}_{N_G(H)\backslash G/N_G(G_{v_{r,s}})}(G_{v_{r,s}})$ 
corresponding to $r$-th subgroup $G_{v_r}\leq G$ up to conjugacy 
%%%%%%%%%%%%%%%%%
via the function 
{\tt ConjugacyClassesSubgroupsNGHOrbitRep(ConjugacyClassesSubgroups($G$),$H$)}.
Then we can get the minimal elements of the $G_{v_{r,s}}$'s 
with ${\rm Ker}\, \psi_1/\varphi_1({\rm Ker}\, \psi_2)=0$ 
via the function {\tt MinConjugacyClassesSubgroups($l$)}. 
Finally, we get a necessary and sufficient condition for 
${\rm Obs}_1(L/K/k)$ $=$ $1$ for each case 
%by the case-by-case analysis 
(see Example \ref{ex12-1}).

{\rm (1-2)} 
$m=2,3,4,16,18,33,43,70,71,130,132,133,176,179,234$ ($15$ cases). 
%%%
Applying {\tt FirstObstructionN($G$)} and {\tt FirstObstructionDnr($G$)}, 
we have ${\rm Obs}_1(L/K/k)=1$. 
Hence we just apply Theorem \ref{thDra89}. 
We have the Schur multiplier $M(G)\simeq\bZ/2\bZ$ 
for $G=12Tm$ $(m=2,3,4,16,18,33,43,70,71,130,132,133,176,179,234)$. 
We obtain a Schur cover 
$1\to M(G)\simeq\bZ/2\bZ\to \widetilde{G}\xrightarrow{\pi} G\to 1$ 
and ${\rm Obs}(K/k)={\rm Obs}_1(\widetilde{L}/K/k)$. 
By Theorem \ref{thmain1}, we have 
${\rm Ker}\, \widetilde{\psi}_1/\widetilde{\varphi}_1({\rm Ker}\, \widetilde{\psi}_2^{\rm nr})\simeq \bZ/2\bZ$. 
By applying the functions 
{\tt FirstObstructionDr($\widetilde{G},\widetilde{G}_{v_{r,s}},\widetilde{H}$)} and {\tt MinConjugacyClassesSubgroups($l$)} 
as in the case (1-1), 
we can get the minimal elements of the $\widetilde{G}_{v_{r,s}}$'s 
with ${\rm Ker}\, \widetilde{\psi}_1/\widetilde{\varphi}_1({\rm Ker}\, \widetilde{\psi}_2)=0$. 
Then we get a necessary and sufficient condition for 
${\rm Obs}_1(\widetilde{L}/K/k)$ $=$ $1$ for each case 
%by the case-by-case analysis 
(see Example \ref{ex12-2}).

{\rm (1-3)} 
$m=20,40,117,168,171,194,261,280$ ($8$ cases).  
%%%
Applying the functions 
{\tt FirstObstructionN($G$)} and {\tt FirstObstructionDnr($G$)}, 
we have ${\rm Obs}_1(L/K/k)=1$. 
For the $8$ cases 
$G=12Tm$ $(m=20,40,117,168,171,194$, $261,280)$, 
we obtain that the Schur multipliers 
$M(G)\simeq \bZ/6\bZ$, $\bZ/6\bZ$, 
$(\bZ/2\bZ)^{\oplus 3}$, $(\bZ/2\bZ)^{\oplus 3}$, 
$\bZ/2\bZ\oplus (\bZ/3\bZ)^{\oplus 2}$, 
$\bZ/6\bZ$, $(\bZ/2\bZ)^{\oplus 4}$, $(\bZ/2\bZ)^{\oplus 2}$ 
respectively. 
We can take a minimal stem extension 
$\overline{G}=\widetilde{G}/A^\prime$, 
i.e. $\overline{A}\leq Z(\overline{G})\cap [\overline{G},\overline{G}]$, 
of $G$ in the commutative diagram 
\begin{align*}
\begin{CD}
1 @>>> A=M(G) @>>> \widetilde{G}@>\pi>> G@>>> 1\\
  @. @VVV @VVV @|\\
1 @>>> \overline{A}=A/A^\prime @>>> \overline{G}=\widetilde{G}/A^\prime @>\overline{\pi}>> G@>>> 1
\end{CD}
\end{align*}
with $\overline{A}\simeq\bZ/2\bZ$ or $\bZ/3\bZ$
via the function 
{\tt MinimalStemExtensions($G$)[$j$].MinimalStemExtension}. 
Then we apply Theorem \ref{thDP87} instead of Theorem \ref{thDra89}. 
Applying 
{\tt KerResH3Z($G,H$:iterator)}, 
we can get  ${\rm Ker}\{H^3(\overline{G}_j,\bZ)\xrightarrow{\rm res}\oplus_{i=1}^{m^\prime} H^3(G_i,\bZ)\}=0$ 
for exactly one $j$ with $\overline{A}\simeq\bZ/2\bZ$. 
Because 
$\oplus_{i=1}^{m^\prime} \widehat{H}^{-3}(G_i,\bZ)\xrightarrow{\rm cores} 
\widehat{H}^{-3}(\overline{G}_j,\bZ)$ is surjective, 
it follows from Theorem \ref{thDP87} that 
${\rm Obs}(K/k)={\rm Obs}_1(\overline{L}_j/K/k)$. 
We also checked that 
${\rm Ker}\, \overline{\psi}_1/\overline{\varphi}_1({\rm Ker}\, \overline{\psi}_2^{\rm nr})\simeq\bZ/2\bZ$ for $\overline{G}_j$ and 
${\rm Ker}\, \overline{\psi}_1/\overline{\varphi}_1({\rm Ker}\, \overline{\psi}_2^{\rm nr})=0$ for $\overline{G}_{j^\prime}$ $(j^\prime\neq j)$. 
Hence 
${\rm Obs}(K/k)\neq {\rm Obs}_1(\overline{L}_j/K/k)$ 
when $\overline{L}_j/k$ is unramified. 
By applying the functions 
{\tt FirstObstructionDr($\overline{G}_j,\overline{G}_{j,v_{r,s}},\overline{H}_j$)} and {\tt MinConjugacyClassesSubgroups($l$)} 
as in the case (1-1), 
we may find the minimal elements of the $\overline{G}_{j,v_{r,s}}$'s 
with ${\rm Ker}\, \overline{\psi}_1/\overline{\varphi}_1({\rm Ker}\, \overline{\psi}_2)=0$. 
Then we get a necessary and sufficient condition for 
${\rm Obs}_1(\overline{L}_j/K/k)$ $=$ $1$ for each case 
%by the case-by-case analysis 
(see Example \ref{ex12-3}). 

%%%%%%%%%%%%%%%%%%%%%%%%%%%%%%%%%%%%%%%%%%%%%%%%%%%%%%%%%%%%%%%%%%%%%
{\rm (2)} Table $2$-$2$: 
$G=12Tm$ $(m=10,37,54,56,59,61,66,77,88,92,93,100,102,144,188,210,214,230,242,255,271)$ $($$21$ {\rm cases}$)$ with $\Sha(T)\leq \bZ/2\bZ$. 
%We assume that $H$ is the stabilizer of one of the letters is $G$. 

We split the $21$ cases into $2$ parts (2-1), (2-2) 
according to the method to prove the assertion: 

{\rm (2-1)} $m=54,56,59,61,66,88,92,93,100,102,144,188,230,255$ 
($14$ cases). 
This case can be proved by the similar manner as in Case {\rm (1-1)} 
(see Example \ref{ex12-4}). 

{\rm (2-2)} 
$m=10,37,77,210,214,242,271$ ($7$ cases). 
This case can be proved by the similar manner as in Case {\rm (1-3)}. 
In particular, for the $7$ cases 
$G=12Tm$ $(m=10,37,77,210,214,242,271)$, 
we obtain that the Schur multipliers 
$M(G)\simeq (\bZ/2\bZ)^{\oplus 3}$, $(\bZ/2\bZ)^{\oplus 3}$, 
$(\bZ/2\bZ)^{\oplus 3}$, $(\bZ/2\bZ)^{\oplus 3}$, 
$(\bZ/2\bZ)^{\oplus 2}\oplus \bZ/3\bZ$, 
$(\bZ/2\bZ)^{\oplus 5}$, $(\bZ/2\bZ)^{\oplus 3}$ 
respectively (see Example \ref{ex12-5}). 

{\rm (3)} Table $2$-$3$: $G=12T31$ with $\Sha(T)\leq \bZ/4\bZ$. 
This case can be proved by the similar manner as in Case {\rm (1-1)} 
(see Example \ref{ex12-6}). 

{\rm (4)} Table $2$-$4$: $G=12T57$ with $\Sha(T)\leq \bZ/4\bZ$. 
%We assume that $H$ is the stabilizer of one of the letters is $G$. 
This case can be proved by the similar manner as in Case {\rm (1-1)} 
(see Example \ref{ex12-7}). 

{\rm (5)} Table $2$-$5$: $G=12T32$ with $\Sha(T)\leq (\bZ/2\bZ)^{\oplus 2}$. 
%We assume that $H$ is the stabilizer of one of the letters is $G$.\\
This case can be proved by the similar manner as in Case {\rm (1-1)} 
(see Example \ref{ex12-0}). 

The last assertion for the cases 
(1) Table $2$-$1$ and (3) Table $2$-$3$ 
follows because the statements can be described in terms of 
the characteristic subgroups of $G$, i.e. 
invariants under the automorphisms of $G$ 
%%%%%%%%%%%%%%%%%%%%%%%%%%%%%%%%%%%%%%%%%%%%%%%%%%%
although for the cases 
(2) Table $2$-$2$, (4) Table $2$-$4$ and (5) Table $2$-$5$, 
the statements are not stable under the action of all 
the automorphisms ${\rm Aut}(G)$ of $G$. 
We can check this via the function {\tt IsInvariantUnderAutG($l$)}.
\qed\\
%%%%%%%%%%%%%%%%%%%%%%%%%%%%%%%%%%%%%%%%%%%%%%%%%%%%%%%%%%%%%%%%%%%%%%

\smallskip
\begin{example}[$G=12Tm$ $(m=7,9,34,47,52,55,64,65,74,75,96,97,122,132,172,174,232,246)$]\label{ex12-1}~\\\vspace*{-2mm}

(1-1-1) $G=12T7\simeq A_4(6)\times C_2$. 
{\small 

}~\\\vspace*{-4mm}

(1-1-2) $G=12T9\simeq S_4$.
{\small 

}~\\\vspace*{-4mm}

(1-1-3) $G=12T34\simeq (S_3)^2\rtimes C_2$. 
{\small 
\begin{verbatim}
gap> Read("HNP.gap");
gap> G:=TransitiveGroup(12,34);
F_36:2(12e)
gap> H:=Stabilizer(G,1);
Group([ (2,10)(3,11)(4,8)(5,9), (4,8)(5,9)(6,10)(7,11) ])
gap> FirstObstructionN(G,H).ker;
[ [ 2 ], [ [ 2 ], [ [ 1 ] ] ] ]
gap> FirstObstructionDnr(G,H).Dnr;
[ [  ], [ [ 2 ], [  ] ] ]
gap> KerResH3Z(G,H:iterator);
[ [  ], [ [ 2 ], [  ] ] ]
gap> HNPtruefalsefn:=x->FirstObstructionDr(G,x,H).Dr[1]=[2];
function( x ) ... end
gap> Gcs:=ConjugacyClassesSubgroups(G);;
gap> Length(Gcs);
26
gap> GcsH:=ConjugacyClassesSubgroupsNGHOrbitRep(Gcs,H);;
gap> GcsHHNPtf:=List(GcsH,x->List(x,HNPtruefalsefn));;
gap> Collected(List(GcsHHNPtf,Set));
[ [ [ true ], 2 ], [ [ false ], 24 ] ]
gap> GcsHNPfalse:=List(Filtered([1..Length(Gcs)],
> x->false in GcsHHNPtf[x]),y->Gcs[y]);;
gap> Length(GcsHNPfalse);
24
gap> GcsHNPtrue:=List(Filtered([1..Length(Gcs)],
> x->true in GcsHHNPtf[x]),y->Gcs[y]);;
gap> Length(GcsHNPtrue);
2
gap> Collected(List(GcsHNPfalse,x->StructureDescription(Representative(x))));
[ [ "(C3 x C3) : C2", 1 ], [ "(C3 x C3) : C4", 1 ], [ "1", 1 ], [ "C2", 3 ], 
  [ "C2 x C2", 2 ], [ "C3", 2 ], [ "C3 x C3", 1 ], [ "C3 x S3", 2 ], 
  [ "C4", 1 ], [ "C6", 2 ], [ "D12", 2 ], [ "S3", 4 ], [ "S3 x S3", 2 ] ]
gap> Collected(List(GcsHNPtrue,x->StructureDescription(Representative(x))));
[ [ "(S3 x S3) : C2", 1 ], [ "D8", 1 ] ]
\end{verbatim}
}~\\\vspace*{-4mm}

(1-1-4) $G=12T47\simeq (C_3)^2\rtimes Q_8$. 
{\small 
\begin{verbatim}
gap> Read("HNP.gap");
gap> G:=TransitiveGroup(12,47);
[(1/3.3^3):2]E(4)_4
gap> H:=Stabilizer(G,1);
Group([ (2,6,10)(3,7,11)(4,8,12), (3,11)(4,8)(5,9)(6,10) ])
gap> FirstObstructionN(G,H).ker;
[ [ 2 ], [ [ 2 ], [ [ 1 ] ] ] ]
gap> FirstObstructionDnr(G,H).Dnr;
[ [  ], [ [ 2 ], [  ] ] ]
gap> KerResH3Z(G,H:iterator);
[ [  ], [ [ 3 ], [  ] ] ]
gap> HNPtruefalsefn:=x->FirstObstructionDr(G,x,H).Dr[1]=[2];
function( x ) ... end
gap> Gcs:=ConjugacyClassesSubgroups(G);;
gap> Length(Gcs);
14
gap> GcsH:=ConjugacyClassesSubgroupsNGHOrbitRep(Gcs,H);;
gap> GcsHHNPtf:=List(GcsH,x->List(x,HNPtruefalsefn));;
gap> Collected(List(GcsHHNPtf,Set));
[ [ [ true ], 2 ], [ [ false ], 12 ] ]
gap> GcsHNPfalse:=List(Filtered([1..Length(Gcs)],
> x->false in GcsHHNPtf[x]),y->Gcs[y]);;
gap> Length(GcsHNPfalse);
12
gap> GcsHNPtrue:=List(Filtered([1..Length(Gcs)],
> x->true in GcsHHNPtf[x]),y->Gcs[y]);;
gap> Length(GcsHNPtrue);
2
gap> Collected(List(GcsHNPfalse,x->StructureDescription(Representative(x))));
[ [ "(C3 x C3) : C2", 1 ], [ "(C3 x C3) : C4", 3 ], [ "1", 1 ], [ "C2", 1 ], 
  [ "C3", 1 ], [ "C3 x C3", 1 ], [ "C4", 3 ], [ "S3", 1 ] ]
gap> Collected(List(GcsHNPtrue,x->StructureDescription(Representative(x))));
[ [ "(C3 x C3) : Q8", 1 ], [ "Q8", 1 ] ]
\end{verbatim}
}~\\\vspace*{-4mm}

(1-1-5) $G=12T52\simeq (A_4\times V_4)\rtimes C_2$. 
{\small 

}~\\\vspace*{-4mm}

(1-1-6) $G=12T55\simeq ((C_4)^2\rtimes C_3)\times C_2$. 
{\small 
\begin{verbatim}
gap> Read("HNP.gap");
gap> G:=TransitiveGroup(12,55);
[1/2.4^3]3
gap> H:=Stabilizer(G,1);
Group([ (3,9)(6,12), (2,8)(5,11), (2,11,8,5)(3,6,9,12) ])
gap> FirstObstructionN(G,H).ker;
[ [ 4 ], [ [ 2, 4 ], [ [ 1, 1 ], [ 0, 2 ] ] ] ]
gap> FirstObstructionDnr(G,H).Dnr;
[ [ 2 ], [ [ 2, 4 ], [ [ 0, 2 ] ] ] ]
gap> KerResH3Z(G,H:iterator);
[ [  ], [ [ 4 ], [  ] ] ]
gap> HNPtruefalsefn:=x->FirstObstructionDr(G,x,H).Dr[1]=[4];
function( x ) ... end
gap> Gcs:=ConjugacyClassesSubgroups(G);;
gap> Length(Gcs);
28
gap> GcsH:=ConjugacyClassesSubgroupsNGHOrbitRep(Gcs,H);;
gap> GcsHHNPtf:=List(GcsH,x->List(x,HNPtruefalsefn));;
gap> Collected(List(GcsHHNPtf,Set));
[ [ [ true ], 5 ], [ [ false ], 23 ] ]
gap> GcsHNPfalse:=List(Filtered([1..Length(Gcs)],
> x->false in GcsHHNPtf[x]),y->Gcs[y]);;
gap> Length(GcsHNPfalse);
23
gap> GcsHNPtrue:=List(Filtered([1..Length(Gcs)],
> x->true in GcsHHNPtf[x]),y->Gcs[y]);;
gap> Length(GcsHNPtrue);
5
gap> Collected(List(GcsHNPfalse,x->StructureDescription(Representative(x))));
[ [ "1", 1 ], [ "A4", 1 ], [ "C2", 3 ], [ "C2 x A4", 1 ], [ "C2 x C2", 3 ], 
  [ "C2 x C2 x C2", 1 ], [ "C3", 1 ], [ "C4", 4 ], [ "C4 x C2", 6 ], 
  [ "C4 x C2 x C2", 1 ], [ "C6", 1 ] ]
gap> Collected(List(GcsHNPtrue,x->StructureDescription(Representative(x))));
[ [ "(C4 x C4) : C3", 1 ], [ "C2 x ((C4 x C4) : C3)", 1 ], [ "C4 x C4", 2 ], 
  [ "C4 x C4 x C2", 1 ] ]
\end{verbatim}
}~\\\vspace*{-4mm}

(1-1-7) $G=12T64\simeq ((C_4)^2\rtimes C_3)\rtimes C_2$. 
{\small 

}~\\\vspace*{-4mm}

(1-1-8) $G=12T65\simeq ((C_4)^2\rtimes C_3)\rtimes C_2$.
{\small 

}~\\\vspace*{-4mm}

(1-1-9) $G=12T74\simeq S_5(12)$. 
{\small 

}~\\\vspace*{-4mm}

(1-1-10) $G=12T75\simeq A_5(6)\times C_2$.
{\small 

}~\\\vspace*{-4mm}

(1-1-11) $G=12T96\simeq (((C_4)^2\rtimes C_3)\rtimes C_2)\times C_2$. 
{\small 

}~\\\vspace*{-4mm}

(1-1-12) $G=12T97\simeq (((C_4)^2\rtimes C_3)\rtimes C_2)\times C_2$. 
{\small 

}~\\\vspace*{-4mm}

(1-1-13) $G=12T122\simeq ((C_3)^2\rtimes Q_8)\rtimes C_3$. 
{\small 
\begin{verbatim}
gap> Read("HNP.gap");
gap> G:=TransitiveGroup(12,122);
[(1/3.3^3):2]A(4)_4
gap> H:=Stabilizer(G,1);
Group([ (2,6,10)(3,7,11)(4,8,12), (2,3,4)(6,7,8)(10,11,12), (3,11)(4,8)(5,9)
(6,10) ])
gap> FirstObstructionN(G,H).ker;
[ [ 2 ], [ [ 6 ], [ [ 3 ] ] ] ]
gap> FirstObstructionDnr(G,H).Dnr;
[ [  ], [ [ 6 ], [  ] ] ]
gap> KerResH3Z(G,H:iterator);
[ [  ], [ [ 3 ], [  ] ] ]
gap> HNPtruefalsefn:=x->FirstObstructionDr(G,x,H).Dr[1]=[2];
function( x ) ... end
gap> Gcs:=ConjugacyClassesSubgroups(G);;
gap> Length(Gcs);
20
gap> GcsH:=ConjugacyClassesSubgroupsNGHOrbitRep(Gcs,H);;
gap> GcsHHNPtf:=List(GcsH,x->List(x,HNPtruefalsefn));;
gap> Collected(List(GcsHHNPtf,Set));
[ [ [ true ], 4 ], [ [ false ], 16 ] ]
gap> GcsHNPfalse:=List(Filtered([1..Length(Gcs)],
> x->false in GcsHHNPtf[x]),y->Gcs[y]);;
gap> Length(GcsHNPfalse);
16
gap> GcsHNPtrue:=List(Filtered([1..Length(Gcs)],
> x->true in GcsHHNPtf[x]),y->Gcs[y]);;
gap> Length(GcsHNPtrue);
4
gap> Collected(List(GcsHNPfalse,x->StructureDescription(Representative(x))));
[ [ "(C3 x C3) : C2", 1 ], [ "(C3 x C3) : C3", 1 ], [ "(C3 x C3) : C4", 1 ], 
  [ "(C3 x C3) : C6", 1 ], [ "1", 1 ], [ "C2", 1 ], [ "C3", 3 ], 
  [ "C3 x C3", 3 ], [ "C3 x S3", 1 ], [ "C4", 1 ], [ "C6", 1 ], [ "S3", 1 ] ]
gap> Collected(List(GcsHNPtrue,x->StructureDescription(Representative(x))));
[ [ "((C3 x C3) : Q8) : C3", 1 ], [ "(C3 x C3) : Q8", 1 ], [ "Q8", 1 ], 
  [ "SL(2,3)", 1 ] ]
\end{verbatim}
}~\\\vspace*{-4mm}

(1-1-14) $G=12T172\simeq (C_3)^4\rtimes D_4$. 
{\small 
\begin{verbatim}
gap> Read("HNP.gap");
gap> G:=TransitiveGroup(12,172);
1/2[3^4:2^2]E(4)
gap> H:=Stabilizer(G,1);
Group([ (4,8,12), (2,10)(5,9)(7,11)(8,12), (3,7,11), (2,6,10) ])
gap> FirstObstructionN(G,H).ker;
[ [ 2 ], [ [ 2 ], [ [ 1 ] ] ] ]
gap> FirstObstructionDnr(G,H).Dnr;
[ [  ], [ [ 2 ], [  ] ] ]
gap> KerResH3Z(G,H:iterator);
[ [  ], [ [ 6 ], [  ] ] ]
gap> HNPtruefalsefn:=x->FirstObstructionDr(G,x,H).Dr[1]=[2];
function( x ) ... end
gap> Gcs:=ConjugacyClassesSubgroups(G);;
gap> Length(Gcs);
277
gap> GcsH:=ConjugacyClassesSubgroupsNGHOrbitRep(Gcs,H);;
gap> GcsHHNPtf:=List(GcsH,x->List(x,HNPtruefalsefn));;
gap> Collected(List(GcsHHNPtf,Set));
[ [ [ true ], 6 ], [ [ false ], 271 ] ]
gap> GcsHNPfalse:=List(Filtered([1..Length(Gcs)],
> x->false in GcsHHNPtf[x]),y->Gcs[y]);;
gap> Length(GcsHNPfalse);
271
gap> GcsHNPtrue:=List(Filtered([1..Length(Gcs)],
> x->true in GcsHHNPtf[x]),y->Gcs[y]);;
gap> Length(GcsHNPtrue);
6
gap> Collected(List(GcsHNPfalse,x->StructureDescription(Representative(x))));
[ [ "((C3 x C3) : C2) x ((C3 x C3) : C2)", 2 ], 
  [ "((C3 x C3) : C2) x S3", 8 ], [ "(C3 x C3 x C3 x C3) : C2", 1 ], 
  [ "(C3 x C3 x C3 x C3) : C4", 1 ], [ "(C3 x C3 x C3) : C2", 14 ], 
  [ "(C3 x C3) : C2", 46 ], [ "(C3 x C3) : C4", 7 ], [ "1", 1 ], [ "C2", 3 ], 
  [ "C2 x ((C3 x C3) : C2)", 2 ], [ "C2 x C2", 2 ], [ "C3", 14 ], 
  [ "C3 x ((C3 x C3) : C2)", 8 ], [ "C3 x C3", 44 ], 
  [ "C3 x C3 x ((C3 x C3) : C2)", 2 ], [ "C3 x C3 x C3", 14 ], 
  [ "C3 x C3 x C3 x C3", 1 ], [ "C3 x C3 x S3", 8 ], [ "C3 x S3", 32 ], 
  [ "C4", 1 ], [ "C6", 8 ], [ "C6 x C3", 2 ], [ "D12", 8 ], [ "S3", 22 ], 
  [ "S3 x S3", 20 ] ]
gap> Collected(List(GcsHNPtrue,x->StructureDescription(Representative(x))));
[ [ "(C3 x C3 x C3 x C3) : D8", 1 ], [ "(S3 x S3) : C2", 4 ], [ "D8", 1 ] ]
\end{verbatim}
}~\\\vspace*{-4mm}

(1-1-15) $G=12T174\simeq (C_3)^4\rtimes Q_8$. 
{\small 
\begin{verbatim}
gap> Read("HNP.gap");
gap> G:=TransitiveGroup(12,174);
[3^4:2]E(4)_4
gap> H:=Stabilizer(G,1);
Group([ (4,8,12), (2,10)(5,9)(7,11)(8,12), (3,7,11), (2,6,10) ])
gap> FirstObstructionN(G,H).ker;
[ [ 2 ], [ [ 2 ], [ [ 1 ] ] ] ]
gap> FirstObstructionDnr(G,H).Dnr;
[ [  ], [ [ 2 ], [  ] ] ]
gap> KerResH3Z(G,H:iterator);
[ [  ], [ [ 3, 3, 3 ], [  ] ] ]
gap> HNPtruefalsefn:=x->FirstObstructionDr(G,x,H).Dr[1]=[2];
function( x ) ... end
gap> Gcs:=ConjugacyClassesSubgroups(G);;
gap> Length(Gcs);
157
gap> GcsH:=ConjugacyClassesSubgroupsNGHOrbitRep(Gcs,H);;
gap> GcsHHNPtf:=List(GcsH,x->List(x,HNPtruefalsefn));;
gap> Collected(List(GcsHHNPtf,Set));
[ [ [ true ], 6 ], [ [ false ], 151 ] ]
gap> GcsHNPfalse:=List(Filtered([1..Length(Gcs)],
> x->false in GcsHHNPtf[x]),y->Gcs[y]);;
gap> Length(GcsHNPfalse);
151
gap> GcsHNPtrue:=List(Filtered([1..Length(Gcs)],
> x->true in GcsHHNPtf[x]),y->Gcs[y]);;
gap> Length(GcsHNPtrue);
6
gap> Collected(List(GcsHNPfalse,x->StructureDescription(Representative(x))));
[ [ "(C3 x C3 x C3 x C3) : C2", 1 ], [ "(C3 x C3 x C3 x C3) : C4", 3 ], 
  [ "(C3 x C3 x C3) : C2", 10 ], [ "(C3 x C3) : C2", 40 ], 
  [ "(C3 x C3) : C4", 21 ], [ "1", 1 ], [ "C2", 1 ], [ "C3", 10 ], 
  [ "C3 x C3", 40 ], [ "C3 x C3 x C3", 10 ], [ "C3 x C3 x C3 x C3", 1 ], 
  [ "C4", 3 ], [ "S3", 10 ] ]
gap> Collected(List(GcsHNPtrue,x->StructureDescription(Representative(x))));
[ [ "(C3 x C3 x C3 x C3) : Q8", 1 ], [ "(C3 x C3) : Q8", 4 ], [ "Q8", 1 ] ]
\end{verbatim}
}~\\\vspace*{-4mm}

(1-1-16) $G=12T232\simeq ((C_3)^4\rtimes Q_8)\rtimes C_3$. 
{\small 

}~\\\vspace*{-4mm}

(1-1-17) $G=12T246\simeq (C_3)^4\rtimes ((C_2)^3\rtimes C_4)$.
{\small 

}
\end{example}

%%%%%%%%%%%%%%%%%%%%%%%%%%%%%%%%%%%%%%%%%%%%%%%%%%%%%%%%%%%%%%%%%%%%
\smallskip
\begin{example}[$G=12Tm$ $(m=2,3,4,16,18,33,43,70,71,130,132,133,176,179,234)$]\label{ex12-2}~\\\vspace*{-2mm}

(1-2-1) $G=12T2\simeq C_6\times C_2$.
{\small 
\begin{verbatim}
gap> Read("HNP.gap");
gap> G:=TransitiveGroup(12,2);
E(4)[x]C(3)=6x2
gap> H:=Stabilizer(G,1);
Group(())
gap> FirstObstructionN(G,H).ker;
[ [  ], [ [  ], [  ] ] ]
gap> GroupCohomology(G,3); # H^3(G,Z)=M(G): Schur multiplier of G
[ 2 ]
gap> ScG:=SchurCoverG(G);
rec( SchurCover := Group([ (1,2,3)(4,5,6)(7,8,9)(10,11,12), (1,4)(2,5)(3,6)
  (7,10)(8,11)(9,12), (4,7)(5,8)(6,9) ]), Tid := [ 12, 14 ], 
  epi := [ (1,2,3)(4,5,6)(7,8,9)(10,11,12), (1,4)(2,5)(3,6)(7,10)(8,11)(9,12),
      (4,7)(5,8)(6,9) ] -> [ (1,5,9)(2,6,10)(3,7,11)(4,8,12), 
      (1,7)(2,8)(3,9)(4,10)(5,11)(6,12), (1,10)(2,5)(3,12)(4,7)(6,9)(8,11) ] )
gap> tG:=ScG.SchurCover;
Group([ (1,2,3)(4,5,6)(7,8,9)(10,11,12), (1,4)(2,5)(3,6)(7,10)(8,11)
(9,12), (4,7)(5,8)(6,9) ])
gap> StructureDescription(tG);
"C3 x D8"
gap> tH:=PreImage(ScG.epi,H);
Group([ (1,10)(2,11)(3,12)(4,7)(5,8)(6,9) ])
gap> FirstObstructionN(tG,tH).ker;
[ [ 2 ], [ [ 2 ], [ [ 1 ] ] ] ]
gap> FirstObstructionDnr(tG,tH).Dnr;
[ [  ], [ [ 2 ], [  ] ] ]
gap> HNPtruefalsefn:=x->FirstObstructionDr(tG,PreImage(ScG.epi,x),tH).Dr[1]=[2];
function( x ) ... end
gap> Gcs:=ConjugacyClassesSubgroups(G);;
gap> Length(Gcs);
10
gap> GcsH:=ConjugacyClassesSubgroupsNGHOrbitRep(Gcs,H);;
gap> GcsHHNPtf:=List(GcsH,x->List(x,HNPtruefalsefn));;
gap> Collected(List(GcsHHNPtf,Set));
[ [ [ true ], 2 ], [ [ false ], 8 ] ]
gap> GcsHNPfalse:=List(Filtered([1..Length(Gcs)],
> x->false in GcsHHNPtf[x]),y->Gcs[y]);;
gap> Length(GcsHNPfalse);
8
gap> GcsHNPtrue:=List(Filtered([1..Length(Gcs)],
> x->true in GcsHHNPtf[x]),y->Gcs[y]);;
gap> Length(GcsHNPtrue);
2
gap> Collected(List(GcsHNPfalse,x->StructureDescription(Representative(x))));
[ [ "1", 1 ], [ "C2", 3 ], [ "C3", 1 ], [ "C6", 3 ] ]
gap> Collected(List(GcsHNPtrue,x->StructureDescription(Representative(x))));
[ [ "C2 x C2", 1 ], [ "C6 x C2", 1 ] ]
\end{verbatim}
}~\\\vspace*{-4mm}

(1-2-2) $G=12T3\simeq D_6$. 
{\small 
\begin{verbatim}
gap> Read("HNP.gap");
gap> G:=TransitiveGroup(12,3);
D_6(6)[x]2=1/2[3:2]E(4)
gap> H:=Stabilizer(G,1);
Group(())
gap> FirstObstructionN(G,H).ker;
[ [  ], [ [  ], [  ] ] ]
gap> GroupCohomology(G,3); # H^3(G,Z)=M(G): Schur multiplier of G
[ 2 ]
gap> ScG:=SchurCoverG(G);
rec( SchurCover := Group([ (5,6), (2,3)(4,5)(6,7), (1,2,3) ]), 
  epi := [ (5,6), (2,3)(4,5)(6,7), (1,2,3) ] -> 
    [ (1,12)(2,3)(4,5)(6,7)(8,9)(10,11), (1,11)(2,8)(3,9)(4,6)(5,7)(10,12), 
      (1,5,9)(2,6,10)(3,7,11)(4,8,12) ] )
gap> tG:=ScG.SchurCover;
Group([ (5,6), (2,3)(4,5)(6,7), (1,2,3) ])
gap> tH:=PreImage(ScG.epi,H);
Group([ (4,7)(5,6) ])
gap> FirstObstructionN(tG,tH).ker;
[ [ 2 ], [ [ 2 ], [ [ 1 ] ] ] ]
gap> FirstObstructionDnr(tG,tH).Dnr;
[ [  ], [ [ 2 ], [  ] ] ]
gap> HNPtruefalsefn:=x->FirstObstructionDr(tG,PreImage(ScG.epi,x),tH).Dr[1]=[2];
function( x ) ... end
gap> Gcs:=ConjugacyClassesSubgroups(G);;
gap> Length(Gcs);
10
gap> GcsH:=ConjugacyClassesSubgroupsNGHOrbitRep(Gcs,H);;
gap> GcsHHNPtf:=List(GcsH,x->List(x,HNPtruefalsefn));;
gap> Collected(List(GcsHHNPtf,Set));
[ [ [ true ], 2 ], [ [ false ], 8 ] ]
gap> GcsHNPfalse:=List(Filtered([1..Length(Gcs)],
> x->false in GcsHHNPtf[x]),y->Gcs[y]);;
gap> Length(GcsHNPfalse);
8
gap> GcsHNPtrue:=List(Filtered([1..Length(Gcs)],
> x->true in GcsHHNPtf[x]),y->Gcs[y]);;
gap> Length(GcsHNPtrue);
2
gap> Collected(List(GcsHNPfalse,x->StructureDescription(Representative(x))));
[ [ "1", 1 ], [ "C2", 3 ], [ "C3", 1 ], [ "C6", 1 ], [ "S3", 2 ] ]
gap> Collected(List(GcsHNPtrue,x->StructureDescription(Representative(x))));
[ [ "C2 x C2", 1 ], [ "D12", 1 ] ]
\end{verbatim}
}~\\\vspace*{-4mm}

(1-2-3) $G=12T4\simeq A_4(12)$. 
{\small 
\begin{verbatim}
gap> Read("HNP.gap");
gap> G:=TransitiveGroup(12,4);
A_4(12)
gap> H:=Stabilizer(G,1);
Group(())
gap> FirstObstructionN(G,H).ker;
[ [  ], [ [  ], [  ] ] ]
gap> GroupCohomology(G,3); # H^3(G,Z)=M(G): Schur multiplier of G
[ 2 ]
gap> ScG:=SchurCoverG(G);
rec( SchurCover := Group([ (2,4,5)(3,7,8), (1,2,6,3)(4,5,7,8) ]), 
  Tid := [ 8, 12 ], epi := [ (2,4,5)(3,7,8), (1,2,6,3)(4,5,7,8) ] -> 
    [ (1,9,5)(2,4,3)(6,8,7)(10,12,11), (1,4)(2,11)(3,6)(5,8)(7,10)(9,12) ] )
gap> tG:=ScG.SchurCover;
Group([ (2,4,5)(3,7,8), (1,2,6,3)(4,5,7,8) ])
gap> StructureDescription(tG);
"SL(2,3)"
gap> tH:=PreImage(ScG.epi,H);
Group([ (1,6)(2,3)(4,7)(5,8) ])
gap> FirstObstructionN(tG,tH).ker;
[ [ 2 ], [ [ 2 ], [ [ 1 ] ] ] ]
gap> FirstObstructionDnr(tG,tH).Dnr;
[ [  ], [ [ 2 ], [  ] ] ]
gap> HNPtruefalsefn:=x->FirstObstructionDr(tG,PreImage(ScG.epi,x),tH).Dr[1]=[2];
function( x ) ... end
gap> Gcs:=ConjugacyClassesSubgroups(G);;
gap> Length(Gcs);
5
gap> GcsH:=ConjugacyClassesSubgroupsNGHOrbitRep(Gcs,H);;
gap> GcsHHNPtf:=List(GcsH,x->List(x,HNPtruefalsefn));;
gap> Collected(List(GcsHHNPtf,Set));
[ [ [ true ], 2 ], [ [ false ], 3 ] ]
gap> GcsHNPfalse:=List(Filtered([1..Length(Gcs)],
> x->false in GcsHHNPtf[x]),y->Gcs[y]);;
gap> Length(GcsHNPfalse);
3
gap> GcsHNPtrue:=List(Filtered([1..Length(Gcs)],
> x->true in GcsHHNPtf[x]),y->Gcs[y]);;
gap> Length(GcsHNPtrue);
2
gap> Collected(List(GcsHNPfalse,x->StructureDescription(Representative(x))));
[ [ "1", 1 ], [ "C2", 1 ], [ "C3", 1 ] ]
gap> Collected(List(GcsHNPtrue,x->StructureDescription(Representative(x))));
[ [ "A4", 1 ], [ "C2 x C2", 1 ] ]
\end{verbatim}
}~\\\vspace*{-4mm}

(1-2-4) $G=12T16\simeq (S_3)^2$.
{\small 
\begin{verbatim}
gap> Read("HNP.gap");
gap> G:=TransitiveGroup(12,16);
[3^2]E(4)
gap> H:=Stabilizer(G,1);
Group([ (2,10,6)(4,8,12) ])
gap> FirstObstructionN(G,H).ker;
[ [ 3 ], [ [ 3 ], [ [ 1 ] ] ] ]
gap> FirstObstructionDnr(G,H).Dnr;
[ [ 3 ], [ [ 3 ], [ [ 1 ] ] ] ]
gap> GroupCohomology(G,3); # H^3(G,Z)=M(G): Schur multiplier of G
[ 2 ]
gap> ScG:=SchurCoverG(G);
rec( SchurCover := Group([ (2,3)(7,8)(9,10), (1,2)(3,4)(5,6)(7,9)
  (8,10), (6,7,8) ]), 
  epi := [ (2,3)(7,8)(9,10), (1,2)(3,4)(5,6)(7,9)(8,10), (6,7,8) ] -> 
    [ (1,7)(2,8)(3,9)(4,10)(5,11)(6,12), (1,10)(2,5)(3,12)(4,7)(6,9)(8,11), 
      (2,10,6)(4,8,12) ] )
gap> tG:=ScG.SchurCover;
Group([ (2,3)(7,8)(9,10), (1,2)(3,4)(5,6)(7,9)(8,10), (6,7,8) ])
gap> StructureDescription(tG);
"(C6 x S3) : C2"
gap> tH:=PreImage(ScG.epi,H);
Group([ (6,7,8), (1,4)(2,3) ])
gap> FirstObstructionN(tG,tH).ker;
[ [ 6 ], [ [ 6 ], [ [ 1 ] ] ] ]
gap> FirstObstructionDnr(tG,tH).Dnr;
[ [ 3 ], [ [ 6 ], [ [ 2 ] ] ] ]
gap> HNPtruefalsefn:=
> x->Product(FirstObstructionDr(tG,PreImage(ScG.epi,x),tH).Dr[1]) mod 2=0;
function( x ) ... end
gap> Gcs:=ConjugacyClassesSubgroups(G);;
gap> Length(Gcs);
22
gap> GcsH:=ConjugacyClassesSubgroupsNGHOrbitRep(Gcs,H);;
gap> GcsHHNPtf:=List(GcsH,x->List(x,HNPtruefalsefn));;
gap> Collected(List(GcsHHNPtf,Set));
[ [ [ true ], 4 ], [ [ false ], 18 ] ]
gap> GcsHNPfalse:=List(Filtered([1..Length(Gcs)],
> x->false in GcsHHNPtf[x]),y->Gcs[y]);;
gap> Length(GcsHNPfalse);
18
gap> GcsHNPtrue:=List(Filtered([1..Length(Gcs)],
> x->true in GcsHHNPtf[x]),y->Gcs[y]);;
gap> Length(GcsHNPtrue);
4
gap> Collected(List(GcsHNPfalse,x->StructureDescription(Representative(x))));
[ [ "(C3 x C3) : C2", 1 ], [ "1", 1 ], [ "C2", 3 ], [ "C3", 3 ], 
  [ "C3 x C3", 1 ], [ "C3 x S3", 2 ], [ "C6", 2 ], [ "S3", 5 ] ]
gap> Collected(List(GcsHNPtrue,x->StructureDescription(Representative(x))));
[ [ "C2 x C2", 1 ], [ "D12", 2 ], [ "S3 x S3", 1 ] ]
\end{verbatim}
}~\\\vspace*{-4mm}

(1-2-5) $G=12T18\simeq S_3\times C_6$.
{\small 

}~\\\vspace*{-4mm}

(1-2-6) $G=12T33\simeq A_5(12)$. 
{\small 

}~\\\vspace*{-4mm}

(1-2-7) $G=12T43\simeq A_4(4)\times S_3$.
{\small 

}~\\\vspace*{-4mm}

(1-2-8) $G=12T70\simeq (S_3)^2\times C_3$. 
{\small 

}~\\\vspace*{-4mm}

(1-2-9) $G=12T71\simeq (C_3)^3\rtimes V_4$. 
{\small 

}~\\\vspace*{-4mm}

(1-2-10) $G=12T130\simeq (C_3)^4\rtimes V_4\simeq C_3\wr V_4$. 
{\small 

}~\\\vspace*{-4mm}

(1-2-11) $G=12T132\simeq ((C_3)^3\rtimes V_4)\rtimes C_3$. 
{\small 

}~\\\vspace*{-4mm}

(1-2-12) $G=12T133\simeq (C_3)^3\rtimes A_4(4)$. 
{\small 

}~\\\vspace*{-4mm}

(1-2-13) $G=12T176\simeq ((C_3)^3\rtimes C_2)\rtimes A_4(4)$. 
{\small 

}~\\\vspace*{-4mm}

(1-2-14) $G=12T179\simeq \PSL_2(\bF_{11})$. 
{\small 

}~\\\vspace*{-4mm}

(1-2-15) $G=12T234\simeq ((C_3)^4\rtimes C_2)\rtimes A_4(4)$. 
{\small 

}
\end{example}

%%%%%%%%%%%%%%%%%%%%%%%%%%%%%%%%%%%%%%%%%%%%%%%%%%%%%%%%%%%%%%%%%%%%
\smallskip
\begin{example}[$G=12Tm$ $(m=20,40,117,168,171,194,261,280)$]\label{ex12-3}~\\\vspace*{-2mm}

(1-3-1) $G=12T20\simeq A_4(4)\times C_3$.
{\small 

}~\\\vspace*{-4mm}

(2-3-2) $G=12T40\simeq ((C_3)^2\rtimes C_4)\times C_2$. 
{\small 
\begin{verbatim}
gap> Read("HNP.gap");
gap> G:=TransitiveGroup(12,40);
F_36(6)[x]2
gap> H:=Stabilizer(G,1);
Group([ (2,10)(3,11)(4,8)(5,9), (2,6,10)(3,7,11) ])
gap> FirstObstructionN(G,H).ker;
[ [  ], [ [ 2 ], [  ] ] ]
gap> GroupCohomology(G,3); # H^3(G,Z)=M(G): Schur multiplier of G
[ 2, 3 ]
gap> cGs:=MinimalStemExtensions(G);;
gap> for cG in cGs do
> bG:=cG.MinimalStemExtension;
> bH:=PreImage(cG.epi,H);
> Print(FirstObstructionN(bG,bH).ker[1]);
> Print(FirstObstructionDnr(bG,bH).Dnr[1]);
> Print("\n");
> od;
[ 3 ][ 3 ]
[ 2 ][  ]
gap> cG:=cGs[2];;
gap> bG:=cG.MinimalStemExtension;
<permutation group of size 144 with 6 generators>
gap> bH:=PreImage(cG.epi,H);
<permutation group of size 12 with 3 generators>
gap> KerResH3Z(bG,bH:iterator);
[ [  ], [ [ 2, 6 ], [  ] ] ]
gap> FirstObstructionN(bG,bH).ker;
[ [ 2 ], [ [ 2, 2 ], [ [ 1, 0 ] ] ] ]
gap> FirstObstructionDnr(bG,bH).Dnr;
[ [  ], [ [ 2, 2 ], [  ] ] ]
gap> HNPtruefalsefn:=x->FirstObstructionDr(bG,PreImage(cG.epi,x),bH).Dr[1]=[2];
function( x ) ... end
gap> Gcs:=ConjugacyClassesSubgroups(G);;
gap> Length(Gcs);
26
gap> GcsH:=ConjugacyClassesSubgroupsNGHOrbitRep(Gcs,H);;
gap> GcsHHNPtf:=List(GcsH,x->List(x,HNPtruefalsefn));;
gap> Collected(List(GcsHHNPtf,Set));
[ [ [ true ], 2 ], [ [ false ], 24 ] ]
gap> GcsHNPfalse:=List(Filtered([1..Length(Gcs)],
> x->false in GcsHHNPtf[x]),y->Gcs[y]);;
gap> Length(GcsHNPfalse);
24
gap> GcsHNPtrue:=List(Filtered([1..Length(Gcs)],
> x->true in GcsHHNPtf[x]),y->Gcs[y]);;
gap> Length(GcsHNPtrue);
2
gap> Collected(List(GcsHNPfalse,x->StructureDescription(Representative(x))));
[ [ "(C3 x C3) : C2", 2 ], [ "(C3 x C3) : C4", 2 ], [ "1", 1 ], [ "C2", 3 ], 
  [ "C2 x ((C3 x C3) : C2)", 1 ], [ "C2 x C2", 1 ], [ "C3", 2 ], 
  [ "C3 x C3", 1 ], [ "C4", 2 ], [ "C6", 2 ], [ "C6 x C3", 1 ], [ "D12", 2 ], 
  [ "S3", 4 ] ]
gap> Collected(List(GcsHNPtrue,x->StructureDescription(Representative(x))));
[ [ "C2 x ((C3 x C3) : C4)", 1 ], [ "C4 x C2", 1 ] ]
\end{verbatim}
}~\\\vspace*{-4mm}

(1-3-3) $G=12T117\simeq (S_3)^3$. 
{\small 

}~\\\vspace*{-4mm}

(1-3-4) $G=12T168\simeq (C_3)^4\rtimes (C_2)^3$. 
{\small 

}~\\\vspace*{-4mm}

(1-3-5) $G=12T171\simeq (C_3)^4\rtimes (C_4\times C_2)$. 
{\small 

}~\\\vspace*{-4mm}

(1-3-6) $G=12T194\simeq (C_3)^4\rtimes A_4(4)\simeq C_3\wr A_4(4)$. 
{\small 

}~\\\vspace*{-4mm}

(1-3-7) $G=12T261\simeq (S_3)^4\rtimes V_4\simeq S_3\wr V_4$. 
{\small 

}~\\\vspace*{-4mm}

(1-3-8) $G=12T280\simeq (S_3)^4\rtimes A_4(4)\simeq S_3\wr A_4(4)$. 
{\small 

}
\end{example}

%%%%%%%%%%%%%%%%%%%%%%%%%%%%%%%%%%%%%%%%%%%%%%%%%%%%%%%%%%%%%%%%

\smallskip
\begin{example}[$G=12Tm$ $(m=54,56,59,61,66,88,92,93,100,102,144,188,230,255)$]\label{ex12-4}~\\\vspace*{-2mm}

(2-1-1) $G=12T54\simeq (S_4\times C_2)\rtimes C_2$. 
{\small 

}~\\\vspace*{-4mm}

(2-1-2) $G=12T56\simeq ((C_2)^4\rtimes C_3)\times C_2$. 
{\small 

}~\\\vspace*{-4mm}

(2-1-3) $G=12T59\simeq (C_2)^3\rtimes A_4(6)$. 
{\small 

}~\\\vspace*{-4mm}

(2-1-4) $G=12T61\simeq ((C_4)^2\rtimes C_2)\rtimes C_3$. 
{\small 

}~\\\vspace*{-4mm}

(2-1-5) $G=12T66\simeq ((C_2)^4\rtimes C_3)\rtimes C_2$. 
{\small 

}~\\\vspace*{-4mm}

(2-1-6) $G=12T88\simeq (C_2)^4\rtimes A_4(6)$. 
{\small 

}~\\\vspace*{-4mm}

(2-1-7) $G=12T92\simeq (((C_4)^2\rtimes C_2)\rtimes C_3)\times C_2$. 
{\small 

}~\\\vspace*{-4mm}

(2-1-8) $G=12T93\simeq (((C_4\rtimes C_4)\times C_2)\rtimes C_2)\rtimes C_3$. 
{\small 

}~\\\vspace*{-4mm}

(2-1-9) $G=12T100\simeq (((C_2)^4\rtimes C_3)\rtimes C_2)\times C_2$. 
{\small 

}~\\\vspace*{-4mm}

(2-1-10) $G=12T102\simeq ((C_2)^4\rtimes C_3)\rtimes C_4$. 
{\small 

}~\\\vspace*{-4mm}

(2-1-11) $G=12T144\simeq (C_2)^5\rtimes A_4(6)$. 
{\small 

}~\\\vspace*{-4mm}

(2-1-12) $G=12T188\simeq (C_2)^6\rtimes A_4(6)\simeq C_2\wr A_4(6)$. 
{\small 

}~\\\vspace*{-4mm}

(2-1-13) $G=12T230\simeq (C_2)^5\rtimes A_5(6)$. 
{\small 

}~\\\vspace*{-4mm}

(2-1-14) $G=12T255\simeq (C_2)^6\rtimes A_5(6)\simeq C_2\wr A_5(6)$. 
{\small 

}
\end{example}

%%%%%%%%%%%%%%%%%%%%%%%%%%%%%%%%%%%%%%%%%%%%%%%%%%%%%%%%%%%%%%%%%%%
\begin{example}[$G=12Tm$ $(m=10,37,77,210,214,242,271)$]\label{ex12-5}~\\\vspace*{-2mm}

(2-2-1) $G=12T10\simeq S_3\times V_4$. 
{\small 

}~\\\vspace*{-4mm}

(2-2-2) $G=12T37\simeq (S_3)^2\times C_2$. 
{\small 

}~\\\vspace*{-4mm}

(2-2-3) $G=12T77\simeq (S_3)^2\rtimes V_4$. 
{\small 

}~\\\vspace*{-4mm}

(2-2-4) $G=12T210\simeq (C_3)^4\rtimes (D_4\times C_2)$. 
{\small 

}~\\\vspace*{-4mm}

(2-2-5) $G=12T214\simeq (C_3)^4\rtimes ((C_4\times C_2)\rtimes C_2)$. 
{\small 

}~\\\vspace*{-4mm}

(2-2-6) $G=12T242\simeq (C_3)^4\rtimes ((C_2)^3\rtimes V_4)$. 
{\small 

}~\\\vspace*{-4mm}

(2-2-7) $G=12T271\simeq ((C_3)^4\rtimes (C_2)^3)\rtimes A_4(4)$. 
{\small 

}
\end{example}

%%%%%%%%%%%%%%%%%%%%%%%%%%%%%%%%%%%%%%%%%%%%%%%%%%%%%%%%%%%%%%%%%%%
\begin{example}[{$G=12T31$}]\label{ex12-6}~\\\vspace*{-2mm}

(3) $G=12T31\simeq (C_4)^2\rtimes C_3$. 
{\small 
\begin{verbatim}
gap> Read("HNP.gap");
gap> G:=TransitiveGroup(12,31);
[4^2]3
gap> H:=Stabilizer(G,1);
Group([ (2,11,8,5)(3,6,9,12) ])
gap> FirstObstructionN(G,H).ker;
[ [ 4 ], [ [ 4 ], [ [ 1 ] ] ] ]
gap> FirstObstructionDnr(G,H).Dnr;
[ [  ], [ [ 4 ], [  ] ] ]
gap> KerResH3Z(G,H:iterator);
[ [  ], [ [ 4 ], [  ] ] ]
gap> HNPtruefalsefn:=x->Product(FirstObstructionDr(G,x,H).Dr[1]);
function( x ) ... end
gap> Gcs:=ConjugacyClassesSubgroups(G);;
gap> Length(Gcs);
10
gap> GcsH:=ConjugacyClassesSubgroupsNGHOrbitRep(Gcs,H);;
gap> GcsHHNPtf:=List(GcsH,x->List(x,HNPtruefalsefn));;
gap> Collected(List(GcsHHNPtf,Set));
[ [ [ 1 ], 5 ], [ [ 2 ], 3 ], [ [ 4 ], 2 ] ]
gap> GcsHNPfalse1:=List(Filtered([1..Length(Gcs)],
> x->1 in GcsHHNPtf[x]),y->Gcs[y]);;
gap> Length(GcsHNPfalse1);
5
gap> GcsHNPfalse2:=List(Filtered([1..Length(Gcs)],
> x->2 in GcsHHNPtf[x]),y->Gcs[y]);;
gap> Length(GcsHNPfalse2);
3
gap> GcsHNPtrue:=List(Filtered([1..Length(Gcs)],
> x->4 in GcsHHNPtf[x]),y->Gcs[y]);;
gap> Length(GcsHNPtrue);
2
gap> Collected(List(GcsHNPfalse1,x->StructureDescription(Representative(x))));
[ [ "1", 1 ], [ "C2", 1 ], [ "C3", 1 ], [ "C4", 2 ] ]
gap> Collected(List(GcsHNPfalse2,x->StructureDescription(Representative(x))));
[ [ "A4", 1 ], [ "C2 x C2", 1 ], [ "C4 x C2", 1 ] ]
gap> Collected(List(GcsHNPtrue,x->StructureDescription(Representative(x))));
[ [ "(C4 x C4) : C3", 1 ], [ "C4 x C4", 1 ] ]
\end{verbatim}
}
\end{example}

\begin{example}[{$G=12T57$}]\label{ex12-7}~\\\vspace*{-2mm}

(4) $G=12T57\simeq ((C_4\times C_2)\rtimes C_4)\rtimes C_3$. 
{\small 

}
\end{example}
%%%%%%%%%%%%%%%%%%%%%%%%%%%%%%%%%%%%%%%%%%%%%%%%%%%%%%%%%%%%%%%%%%%%%%

%%%%%%%%%%%%%%%%%%%%%%%%%%%%%%

\smallskip
\begin{example}[{$G=12T32$}]\label{ex12-0}~\vspace*{-2mm}\\

{\rm (5)} $G=12T32\simeq (C_2)^4\rtimes C_3$. 
{\small 

}
\end{example}
%%%%%%%%%%%%%%%%%%%%%%%%%%%%%%%%%%%%%%%%%%%%%%%%%%%%%%%%%%%%%%%%%%

%%%%%%%%%%%%%%%%%%%%%%%%%%%%%%%%%%%%%%%%%%%%%%%%%%%%%%%%%%%%%%%%%%%%%%%%%%
\section{GAP algorithms}\label{S8}
We give GAP algorithms for computing the total obstruction 
${\rm Obs}(K/k)$ and 
the first obstruction ${\rm Obs}_1(L/K/k)$ 
as in Section \ref{S6}. 
The functions which are provided in this section are available from\\
{\tt https://www.math.kyoto-u.ac.jp/\~{}yamasaki/Algorithm/Norm1ToriHNP/}.\\
~{}\vspace*{-4mm}\\
{\small 

}

%%%%%%%%%%%%%%%%%%%%%%%%%%%%%%%%%%%%%%%%%%%%%%%%%%%%%%%%%%%%%%%%%%%%%%%%


\begin{thebibliography}{KMRT98}
\bibitem[Bar81a]{Bar81a} 
H.-J. Bartels, 
{\it Zur Arithmetik von Konjugationsklassen in algebraischen Gruppen}, 
J. Algebra \textbf{70} (1981) 179--199. 
\bibitem[Bar81b]{Bar81b} 
H.-J. Bartels, 
{\it Zur Arithmetik von Diedergruppenerweiterungen}, 
Math. Ann. \textbf{256} (1981) 465--473.
\bibitem[BT82]{BT82} 
F. R. Beyl, J. Tappe, 
{\it Group extensions, representations, and the Schur multiplicator}, 
Lecture Notes in Mathematics, 958. Springer-Verlag, 
Berlin-New York, 1982. 
%\bibitem[But93]{But93} 
%G. Butler, 
%{\it The transitive groups of degree fourteen and fifteen}, 
%J. Symbolic Comput. \textbf{16} (1993) 413--422.
\bibitem[BM83]{BM83}
G. Butler, J. McKay, 
{\it The transitive groups of degree up to eleven}, 
Comm. Algebra \textbf{11} (1983) 863--911.
\bibitem[CT07]{CT07} 
J.-L. Colliot-Th\'{e}l\`{e}ne, 
{\it Lectures on Linear Algebraic Groups}, 
Beijing lectures, Moning side centre, April 2007, 
https://www.math.u-psud.fr/\~{}colliot/BeijingLectures2Juin07.pdf. 
\bibitem[CTHS05]{CTHS05} 
J.-L. Colliot-Th\'{e}l\`{e}ne, D. Harari, A. N. Skorobogatov, 
{\it Compactification \'equivariante d'un tore 
(d'apr\`es Brylinski et K\"unnemann)}, 
Expo. Math. \textbf{23} (2005) 161--170.
\bibitem[CTS77]{CTS77} 
J.-L. Colliot-Th\'{e}l\`{e}ne, J.-J. Sansuc, 
{\it La R-\'{e}quivalence sur les tores}, 
Ann. Sci. \'{E}cole Norm. Sup. (4) \textbf{10} (1977) 175--229. 
\bibitem[CTS87]{CTS87} 
J.-L. Colliot-Th\'{e}l\`{e}ne, J.-J. Sansuc, 
{\it Principal homogeneous spaces under flasque tori: Applications}, 
J. Algebra \textbf{106} (1987) 148--205. 
%%%%%%%%%%%%%%%%%%%%%%%%%%%%%%%%%%%%%%%%%%%%%%%%%%%%
\bibitem[CK00]{CK00} 
A. Cortella, B. Kunyavskii, 
{\it Rationality problem for generic tori in simple groups}, 
J. Algebra \textbf{225} (2000) 771--793. 
\bibitem[Dra89]{Dra89} 
Yu. A. Drakokhrust, 
{\it On the complete obstruction to the Hasse principle}, (Russian)
Dokl. Akad. Nauk BSSR \textbf{30} (1986) 5--8; 
translation in Amer. Math. Soc. Transl. (2) \textbf{143} (1989) 29--34. 
\bibitem[DP87]{DP87} 
Yu. A. Drakokhrust, V. P. Platonov, 
{\it The Hasse norm principle for algebraic number fields}, (Russian)
Izv. Akad. Nauk SSSR Ser. Mat. \textbf{50} (1986) 946--968; 
translation in Math. USSR-Izv. \textbf{29} (1987) 299--322. 
\bibitem[End11]{End11} 
S. Endo, 
{\it The rationality problem for norm one tori}, 
Nagoya Math. J. \textbf{202} (2011) 83--106. 
%\bibitem[EK17]{EK17} 
%S. Endo, M. Kang, 
%{\it Function fields of algebraic tori revisited}, 
%Asian J. Math. \textbf{21} (2017) 197--224.
\bibitem[EM73]{EM73} 
S. Endo, T. Miyata, 
{\it Invariants of finite abelian groups}, 
J. Math. Soc. Japan \textbf{25} (1973) 7--26. 
\bibitem[EM75]{EM75} 
S. Endo, T. Miyata, 
{\it On a classification of the function fields of algebraic tori}, 
Nagoya Math. J. \textbf{56} (1975) 85--104. 
Corrigenda: Nagoya Math. J. \textbf{79} (1980) 187--190. 
\bibitem[Flo]{Flo} 
M. Florence, 
{\it Non rationality of some norm-one tori}, preprint (2006).
%\bibitem[FLN18]{FLN18}
%C. Frei, D. Loughran, R. Newton, 
%{\it The Hasse norm principle for abelian extensions}, 
%Amer. J. Math. \textbf{140} (2018) 1639--1685.
\bibitem[GAP]{GAP} 
The GAP Group, GAP -- Groups, Algorithms, and Programming, 
%Version 4.4.12; 2008. 
%Version 4.8.10; 2018. 
Version 4.9.3; 2018. (http://www.gap-system.org).
\bibitem[Ger77]{Ger77} 
F. Gerth III, 
{\it The Hasse norm principle in metacyclic extensions of number fields}, 
J. London Math. Soc. (2) \textbf{16} (1977) 203--208. 
\bibitem[Ger78]{Ger78} 
F. Gerth III, 
{\it The Hasse norm principle in cyclotomic number fields}, 
J. Reine Angew. Math. \textbf{303/304} (1978) 249--252. 
\bibitem[Gur78a]{Gur78a} 
S. Gurak, {\it On the Hasse norm principle}, 
J. Reine Angew. Math. \textbf{299/300} (1978) 16--27.
\bibitem[Gur78b]{Gur78b}
S. Gurak, 
{\it The Hasse norm principle in non-abelian extensions}, 
J. Reine Angew. Math. \textbf{303/304} (1978) 314--318.
\bibitem[Gur80]{Gur80} 
S. Gurak, 
{\it The Hasse norm principle in a compositum of radical extensions}, 
J. London Math. Soc. (2) \textbf{22} (1980) 385--397.
\bibitem[HAP]{HAP}
G. Ellis, The GAP package HAP, version 1.12.6, available from
http://www.gap-system.org/Packages/hap.html.
\bibitem[HHY20]{HHY20} 
S. Hasegawa, A. Hoshi, A. Yamasaki, 
{\it Rationality problem for norm one tori in small dimensions}, 
Math. Comp. \textbf{89} (2020) 923--940.
\bibitem[Has31]{Has31} 
H. Hasse, 
{\it Beweis eines Satzes und Wiederlegung einer Vermutung \"uber das allgemeine Normenrestsymbol}, Nachrichten von der Gesellschaft der Wissenschaften zu G\"ottingen, Mathematisch-Physikalische Klasse (1931) 64--69. 
\bibitem[Hir64]{Hir64} 
H. Hironaka, {\it Resolution of singularities of an algebraic 
variety over a field of characteristic zero. I, II.}, 
Ann. of Math. (2) \textbf{79} (1964) 109--203; 205--326. 
\bibitem[HKY22]{HKY}
A. Hoshi, K. Kanai, A. Yamasaki, 
{\it Norm one tori and Hasse norm principle}, 
Math. Comp. \textbf{91} (2022) 2431--2458. 
%\bibitem[HKY2]{HKY2}
%A. Hoshi, K. Kanai, A. Yamasaki, 
%{\it Norm one tori and Hasse norm principle, II: Degree $12$ case}, 
%arXiv:2003.08253 (the arXiv version of this paper). 
%%%%%%%%%%%%%%
%\bibitem[HKK14]{HKK14} 
%A. Hoshi, M. Kang, H. Kitayama, 
%{\it Quasi-monomial actions and some 4-dimensional rationality problems}, 
%J. Algebra \textbf{403} (2014) 363--400. 
\bibitem[HY17]{HY17} 
A. Hoshi, A. Yamasaki, 
{\it Rationality problem for algebraic tori}, 
Mem. Amer. Math. Soc. \textbf{248} (2017) no. 1176, v+215 pp. 
\bibitem[HY21]{HY21} 
A. Hoshi, A. Yamasaki, 
{\it Rationality problem for norm one tori}, 
Israel J. Math. \textbf{241} (2021) 849--867. 
\bibitem[H\"{u}r84]{Hur84} 
W. H\"{u}rlimann, 
{\it On algebraic tori of norm type}, Comment. Math. Helv. 
\textbf{59} (1984) 539--549.
\bibitem[Kan12]{Kan12}
M. Kang, {\it Retract rational fields}, 
J. Algebra \textbf{349} (2012) 22--37.
%\bibitem[Mer08]{Mer08} 
%A. Merkurjev, {\it $R$-equivalence on three-dimensional tori 
%and zero-cycles}, 
%Algebra Number Theory \textbf{2} (2008) 69--89.
\bibitem[Kap87]{Kap87} 
G. Karpilovsky, 
{\it The Schur multiplier}, 
London Mathematical Society Monographs. 
New Series, 2. The Clarendon Press, Oxford University Press, New York, 1987.
\bibitem[KMRT98]{KMRT98} 
M.-A. Knus, A. Merkurjev, M. Rost, J.-P. Tignol, 
{\it The book of involutions}, American Mathematical Society Colloquium Publications, 44, American Mathematical Society, Providence, RI, 1998, xxii+593 pp. 
\bibitem[Kun84]{Kun84} 
B. E. Kunyavskii, 
{\it Arithmetic properties of three-dimensional algebraic tori}, (Russian)
Integral lattices and finite linear groups, 
Zap. Nauchn. Sem. Leningrad. Otdel. Mat. Inst. Steklov. (LOMI) \textbf{116} (1982) 102--107, 163; 
translation in J. Soviet Math. \textbf{26} (1984) 1898--1901. 
\bibitem[Kun90]{Kun90}
B. E. Kunyavskii, {\it Three-dimensional algebraic tori},
Selecta Math. Soviet. \textbf{9} (1990) 1--21.
\bibitem[Kun07]{Kun07} 
B. E. Kunyavskii, 
{\it Algebraic tori --- thirty years after}, 
Vestnik Samara State Univ. (2007) 198--214. 
\bibitem[LL00]{LL00} 
N. Lemire, M. Lorenz, 
{\it On certain lattices associated with generic division algebras}, 
J. Group Theory \textbf{3} (2000) 385--405.
\bibitem[LPR06]{LPR06} 
N. Lemire, V. L. Popov, Z. Reichstein, 
{\it Cayley groups}, J. Amer. Math. Soc. \textbf{19} (2006) 921--967. 
\bibitem[LeB95]{LeB95} 
L. Le Bruyn, {\it Generic norm one tori},
Nieuw Arch. Wisk. (4) \textbf{13} (1995) 401--407. 
%\bibitem[Len74]{Len74} 
%H. W. Lenstra, Jr., 
%{\it Rational functions invariant under a finite abelian group}, 
%Invent. Math. \textbf{25} (1974) 299--325.
\bibitem[Lor05]{Lor05} 
M. Lorenz, 
{\it Multiplicative invariant theory}, 
Encyclopaedia Math. Sci., vol. 135, Springer-Verlag, Berlin, 2005. 
\bibitem[Mac20]{Mac20} 
A. Macedo, 
{\it The Hasse norm principle for $A_n$-extensions}, 
J. Number Theory \textbf{211} (2020) 500--512.
\bibitem[Man74]{Man74} 
Yu. I. Manin, 
{\it Cubic forms: algebra, geometry, arithmetic}, 
%Translated from the Russian by M. Hazewinkel. 
%North-Holland Mathematical Library, Vol. 4. 
%North-Holland Publishing Co., Amsterdam-London; 
%American Elsevier Publishing Co., New York, 1974. vii+292 pp.
North-Holland Mathematical Library 4, North-Holland, Amsterdam, 1974.
%\bibitem[Man86]{Man86} Yu. I. Manin, 
%{\it Cubic forms}, 2nd ed., 
%North-Holland Mathematical Library 4, North-Holland, Amsterdam, 1986.
%\bibitem[Maz82]{Maz82} 
%P. Mazet, 
%{\it Sur les multiplicateurs de Schur des groupes de Mathieu}, (French) 
%J. Algebra \textbf{77} (1982) 552--576. 
\bibitem[Ono61]{Ono61} 
T. Ono, 
{\it Arithmetic of algebraic tori}, 
Ann. of Math. (2) \textbf{74} (1961) 101--139. 
\bibitem[Ono63]{Ono63} 
T. Ono, 
{\it On the Tamagawa number of algebraic tori}, 
Ann. of Math. (2) \textbf{78} (1963) 47--73.
%\bibitem[Ono65]{Ono65} 
%T. Ono, 
%{\it On the relative theory of Tamagawa numbers}, 
%Ann. of Math. (2) \textbf{82} (1965) 88--111.
\bibitem[Opo80]{Opo80} 
H. Opolka, 
{\it Zur Aufl\"osung zahlentheoretischer Knoten}, %(German)
Math. Z. \textbf{173} (1980) 95--103. 
\bibitem[Pla82]{Pla82} 
V. P. Platonov, 
{\it Arithmetic theory of algebraic groups}, (Russian)
Uspekhi Mat. Nauk \textbf{37} (1982) 3--54; 
translation in Russian Math. Surveys \textbf{37} (1982) 1--62. 
\bibitem[PD85a]{PD85a} 
V. P. Platonov, Yu. A. Drakokhrust, 
{\it On the Hasse principle for algebraic number fields}, (Russian) 
Dokl. Akad. Nauk SSSR \textbf{281} (1985) 793--797; 
translation in  Soviet Math. Dokl. \textbf{31} (1985) 349--353.
\bibitem[PD85b]{PD85b} 
V. P. Platonov, Yu. A. Drakokhrust, 
{\it The Hasse norm principle for primary extensions of algebraic number fields}, (Russian) Dokl. Akad. Nauk SSSR \textbf{285} (1985) 812--815; 
translation in Soviet Math. Dokl. \textbf{32} (1985) 789--792.
\bibitem[PR94]{PR94} 
V. P. Platonov, A. Rapinchuk, 
{\it Algebraic groups and number theory}, 
Translated from the 1991 Russian original by Rachel Rowen, 
Pure and applied mathematics, 139, Academic Press, 1994. 
\bibitem[Rob96]{Rob96} 
D. J. S. Robinson, 
{\it A course in the theory of groups}, Second edition. 
Graduate Texts in Mathematics, 80. Springer-Verlag, New York, 1996. 
\bibitem[Roy87]{Roy87}
G. F. Royle, 
{\it The transitive groups of degree twelve}, 
J. Symbolic Comput. \textbf{4} (1987) 255--268.
\bibitem[Sal84]{Sal84} 
D. J. Saltman, 
{\it Retract rational fields and cyclic Galois extensions}, 
Israel J. Math. \textbf{47} (1984) 165--215.  
\bibitem[San81]{San81} 
J.-J. Sansuc, 
{\it Groupe de Brauer et arithm\'etique des groupes alg\'ebriques lin\'eaires sur un corps de nombres}, (French) 
J. Reine Angew. Math. \textbf{327} (1981) 12--80.
\bibitem[Swa83]{Swa83} 
R. G. Swan, 
{\it Noether's problem in Galois theory}, 
Emmy Noether in Bryn Mawr (Bryn Mawr, Pa., 1982), 21--40, Springer, 
New York-Berlin, 1983. 
\bibitem[Swa10]{Swa10} 
R. G. Swan, 
{\it The flabby class group of a finite cyclic group}, 
Fourth International Congress of Chinese Mathematicians, 259--269, 
AMS/IP Stud. Adv. Math., 48, Amer. Math. Soc., Providence, RI, 2010. 
\bibitem[Tat67]{Tat67} 
J. Tate, 
{\it Global class field theory}, 
Algebraic Number Theory (Proc. Instructional Conf., Brighton, 1965), 
162--203, Thompson, Washington, D.C., 1967.  
\bibitem[Vos67]{Vos67} 
V. E. Voskresenskii, 
{\it On two-dimensional algebraic tori II}, (Russian) 
Izv. Akad. Nauk SSSR Ser. Mat. \textbf{31} (1967) 711--716; 
translation in Math. USSR-Izv. \textbf{1} (1967) 691--696.
\bibitem[Vos69]{Vos69} 
V. E. Voskresenskii, 
{\it The birational equivalence of linear algebraic groups}, (Russian) 
Dokl. Akad. Nauk SSSR \textbf{188} (1969) 978--981; 
erratum, ibid. 191 1969 nos., 1, 2, 3, vii; 
translation in Soviet Math. Dokl. \textbf{10} (1969) 1212--1215. 
\bibitem[Vos70]{Vos70} 
V. E. Voskresenskii, 
{\it Birational properties of linear algebraic groups}, (Russian) 
Izv. Akad. Nauk SSSR Ser. Mat. \textbf{34} (1970) 3--19; 
translation in Math. USSR-Izv. \textbf{4} (1970) 1--17.
\bibitem[Vos74]{Vos74} 
V. E. Voskresenskii, 
{\it Stable equivalence of algebraic tori}, (Russian) 
Izv. Akad. Nauk SSSR Ser. Mat. \textbf{38} (1974) 3--10; 
translation in Math. USSR-Izv. \textbf{8} (1974) 1--7.
%\bibitem[Vos83]{Vos83} 
%V. E. Voskresenskii, 
%{\it Projective invariant Demazure models}, (Russian)
%Izv. Akad. Nauk SSSR Ser. Mat. \textbf{46} (1982) 195--210, 431; 
%translation in Math USSR-Izv. \textbf{20} (1983) 189--202.
\bibitem[Vos88]{Vos88} 
V. E. Voskresenskii, 
{\it Maximal tori without affect in semisimple algebraic groups}, (Russian) 
Mat. Zametki 44 (1988) 309--318; 
translation in Math. Notes \textbf{44} (1988) 651--655. 
\bibitem[Vos98]{Vos98} 
V. E. Voskresenskii, 
{\it Algebraic groups and their birational invariants}, 
Translated from the Russian manuscript by Boris Kunyavskii, 
Translations of Mathematical Monographs, 179. 
American Mathematical Society, Providence, RI, 1998.
\bibitem[VK84]{VK84} 
V. E. Voskresenskii,  B. E. Kunyavskii, 
{\it Maximal tori in semisimple algebraic groups}, 
Kuibyshev State Inst., Kuibyshev (1984). 
Deposited in VINITI March 5, 1984, No. 1269-84 Dep. 
(Ref. Zh. Mat. (1984), 7A405 Dep.).
\bibitem[Yam12]{Yam12} A. Yamasaki, 
{\it Negative solutions to three-dimensional monomial Noether problem},
J. Algebra \textbf{370} (2012) 46--78.
\end{thebibliography}
\end{document}